\documentclass[12pt,amssymb]{article}
\usepackage{latexsym}

\title {Countable imaginary simple unidimensional theories}
\author {Ziv Shami\\Tel Aviv University}

\newtheorem {theorem}{Theorem}[section]
\newtheorem {lemma}[theorem]{Lemma}
\newtheorem {definition}[theorem]{Definition}

\newtheorem {fact}[theorem]{Fact}

\newtheorem {corollary}[theorem]{Corollary}

\newtheorem {remark}[theorem]{Remark}

\newtheorem {proposition}[theorem]{Proposition}

\newtheorem {claim}[theorem]{Claim}

\newtheorem {subclaim}[theorem]{Subclaim}

\newtheorem {example}[theorem]{Example}

\def\proof {\noindent \textbf{Proof:} }
\def\|#1| {\vert #1 \vert}

\newsavebox{\indbin}
\savebox{\indbin}{\begin{picture}(0,0)
\newlength{\gnu}
\settowidth{\gnu}{$\smile$} \setlength{\unitlength}{.5\gnu} \put(-1,-.65){$\smile$}
\put(-.25,.1){$|$}
\end{picture}}
\newcommand{\nonfork}[3]
{\mbox{$\begin{array}{ccc} \mbox{$#1$} & \usebox{\indbin} & \mbox{$#2$} \\
        & \mbox{$#3$} &
\end{array}$}}
\newcommand{\nonforkempty}[2]
{\mbox{$\begin{array}{ccc} \mbox{$#1$} & \usebox{\indbin} & \mbox{$#2$}
\end{array}$}}
\newcommand{\fork}[3]
{\mbox{$\begin{array}{ccc} \mbox{$#1$} & \!\mbox{$\!\!\not\!\:\usebox{\indbin}$} & \mbox{$#2$} \\
        & \mbox{$#3$} &
\end{array}$}}
\newcommand{\forkempty}[2]
{\mbox{$\begin{array}{ccc} \mbox{$#1$} & \!\mbox{$\!\!\not\!\:\usebox{\indbin}$} & \mbox{$#2$}
\end{array}$}}

\newsavebox{\sindbin}
\savebox{\sindbin}{\begin{picture}(0,0)
\newlength{\sgnu}
\settowidth{\sgnu}{$\smile$} \setlength{\unitlength}{.5\sgnu} \put(-1,-.65){$\smile$}
\put(-.25,.1){$|s$}
\end{picture}}
\newcommand{\snonfork}[3]
{\mbox{$\begin{array}{ccc} \mbox{$#1$} & \usebox{\sindbin} & \mbox{$#2$} \\
        & \mbox{$#3$} &
\end{array}$}}
\newcommand{\snonforkempty}[2]
{\mbox{$\begin{array}{ccc} \mbox{$#1$} & \usebox{\sindbin} & \mbox{$#2$}
\end{array}$}}
\newcommand{\sfork}[3]
{\mbox{$\begin{array}{ccc} \mbox{$#1$} & \!\mbox{$\!\!\not\!\:\usebox{\sindbin}$} & \mbox{$#2$} \\
        & \mbox{$#3$} &
\end{array}$}}

\newsavebox{\starindbin}
\savebox{\starindbin}{\begin{picture}(0,0)
\newlength{\stargnu}
\settowidth{\stargnu}{$\smile$} \setlength{\unitlength}{.5\stargnu} \put(-1,-.65){$\smile$}
\put(-.25,.1){$|*$}
\end{picture}}

\newsavebox{\qindbin}
\savebox{\qindbin}{\begin{picture}(0,0)
\newlength{\qgnu}
\settowidth{\qgnu}{$\smile$} \setlength{\unitlength}{.5\qgnu} \put(-1,-.65){$\smile$}
\put(-.25,.1){$|_{qf}$}
\end{picture}}

\def\card #1 {{\vert #1 \vert}}
\def\acl #1 {{acl^{eq}(#1)}}
\def\CC {{\cal C}}

\def\PP {{\cal P}}
\def\UU {{\cal U}}

\def\WW {{\cal W}}

\begin{document}
\maketitle

\begin{abstract}
We prove that a countable simple unidimensional theory that eliminates hyperimaginaries is
supersimple. This solves a problem of Shelah in the more general context of simple theories under
weak assumptions.
\end{abstract}

\section{Introduction}
The notion of a unidimensional theory already appeared, in a different form, in Baldwin-Lachlan
characterization of $\aleph_1$-categorical theories; a countable theory is $\aleph_1$-categorical
iff it is $\omega$-stable and has no Vaughtian pairs (equivalently, $T$ is $\omega$-stable and
unidimensional). Later, Shelah defined a unidimensional theory to be a stable theory $T$ in which
any two $\vert T\vert ^+$-staurated models of the same power are isomorphic, and proved that in the
stable context a theory is unidimensional iff any two non-algebraic types are non-orthogonal. A
problem posed by Shelah was whether any unidimensional stable theory is superstable. This was
answered positively by Hrushovski around 1986 first in the countable case [H0] and then in full
generality [H1]. Taking the right hand side of Shelah characterization of unidimensional stable
theories seems natural for the simple case. Shelah's problem extended to this context seems much
harder. In [S3] it was observed that a small simple unidimensional theory is supersimple. Later,
Pillay [P] gave a positive answer for countable imaginary simple theories with wnfcp (the weak non
finite cover property), building on the arguments in [H0] and using some machinery from [BPV]. Then
using the result on elimination of $\exists^\infty$ in simple unidimensional theories [S1]
completed his proof for countable imaginary low theories [P1].

In this paper we prove the result for any countable imaginary simple theory. One of the key notions
that will take place in this paper is the forking topology. For variables $x$ and set $A$ the
forking topology on $S_x(A)$ is defined as the topology whose basis is the collection of all sets
of the form $\UU=\{a \vert \phi(a,y)$ forks over $A \}$, where $\phi(x,y)\in L(A)$. These
topologies were defined in [S2] (the $\tau^f$-topologies) and are variants of Hrushovski's [H0] and
later Pillay's [P] topologies. The main role of Hrushovski's and Pillay's topologies in their proof
was the ability to express the relation $\Gamma_F(x)$ defined by $\Gamma_F(x)=\exists y
(F(x,y)\wedge \nonforkempty{y}{x})$ as a closed relation for any Stone-closed relation $F(x,y)$.
Indeed, using this and a property of $T$, we call PCFT, that says these topologies are closed under
projections, they proved the existence of an unbounded $\tau^f$-open set of bounded finite
$SU$-rank in any countable imaginary unidimensional stable/low theory. From this, supersimplicity
follows quite easily by showing that the existence of such a set in a simple theory actually
implies there is a definable set of $SU$-rank 1. In [S2] however, the forking topologies played a
different role. It is shown there, in particular, that if $T$ is an imaginary simple unidimensional
theory with PCFT then the existence of an unbounded supersimple $\tau^f$-open set implies the
theory is supersimple (supersimplicity here does not follow easily as before since we don't know
there is a finite bound on the $SU$-rank of all types extending the supersimple $\tau^f$-open set).

The first step of the proof in the current paper is to show that any simple unidimesional theory
has PCFT. Thus, for proving the main result, it will be sufficient to show there exists an
unbounded supersimple $\tau^f$-open set. The existence of such a set is achieved via the
introduction of the dividing line "$T$ is essentially 1-based" which means every type is
coordinatised by essentially 1-based types in the sense of the forking topology. In case $T$ is not
essentially 1-based we prove there is an unbounded $\tau^f$-open set of finite $SU$-rank (possibly
with no finite bound); this is a general dichotomy for countable imaginary simple theories. If $T$
is essentially 1-based, the problem is reduced to the task of finding an unbounded type-definable
$\tau^f$-open set of bounded finite $SU_s$-rank (the foundation rank with respect to forking with
stable formulas). In order to show the existence of such a set, we introduce the notion of a
$\tilde\tau^f_{st}$-set and prove a theorem saying that in any simple theory in which the extension
property is first-order, any minimal unbounded fiber in an unbounded $\tilde\tau^f_{st}$-set is a
type-definable $\tau^f$-open set. Then, we show that the assumption that $T$ is countable,
imaginary, and unidimensional implies there is a minimal unbounded fiber of some
$\tilde\tau^f_{st}$-set that has bounded finite $SU_s$-rank. By the above theorem we conclude that
this fiber is a type-definable $\tau^f$-open set and thus the proof of the main result is complete.

We will assume basic knowledge of simple theories as in [K1],[KP],[HKP]. A good text book on simple
theories that covers much more is [W]. The notations are standard, and throughout the paper we work
in a highly saturated, highly strongly-homogeneous model $\CC$ of a complete first-order theory $T$
in a language $L$. We will often work in $\CC^{eq}$ and will not work with hyperimaginaries unless
otherwise stated.

\section{Preliminaries}
We recall here some definitions and facts relevant for this paper. In this section $T$ will be a
simple theory and we work in $\CC^{eq}$.

\subsection{Interaction}
For the rest of this section let $\PP$ be an $A$-invariant set of small partial types and $p\in
S(A)$. We say that \em $p$ is (almost-) $\PP$-internal \em if there exists a realization $a$ of $p$
and there exists $B\supseteq A$ with $\nonfork{a}{B}{A}$ such that for some tuple $\bar c$ of
realizations of types in $\PP$ that extend to types in $S(B)$ we have $a \in dcl(B,\bar c)$
(respectively, $a \in acl(B,\bar c)$). We say that \em $p$ is analyzable in $\PP$ \em if there
exists a sequence $I=\langle a_i \vert i\leq \alpha \rangle\subseteq dcl(a_\alpha A)$, where
$a_\alpha\models p$, such that $tp(a_i/A\cup \{a_j \vert j<i\})$ is $\PP$-internal for every $i\leq
\alpha$. We say that $p$ is \em foreign $\PP$ \em if for every $B\supseteq A$ and $a\models p$ with
$\nonfork{a}{B}{A}$ and a realization $c$ of a type in $\PP$ that extends to a type in $S(B)$,
$\nonfork{a}{c}{B}$. Also, recall that $p\in S(A)$ is said to be \em orthogonal \em to some $q\in
S(B)$ if for every $C\supseteq A\cup B$, for every $\bar p\in S(C)$, a non-forking extension of $p$
, and every $\bar q\in S(C)$, a non-forking extension of $q$, for every realization $a$ of $\bar p$
and realization $b$ of $\bar q$, $\nonfork{a}{b}{C}$. The above definitions are valid for
hyperimaginaries as well. Note that in the hyperimaginary context we say that $p$ is analyzable in
$\PP$ (by hyperimaginaries) if there exists a sequence $I\subseteq dcl(a_\alpha A)$ as above of
hyperimaginaries. We say that $T$ is \em imaginary (or has elimination of hyperimaginaries) \em if
for every type-definable over $\emptyset$ equivalence relation $E$ on a complete type $q$ (of a
possibly infinite tuple of elements), $E$ is equivalent on $q$ to the intersection of some
definable equivalence relations $E_i\in L$.

\begin{fact}\label{internal}
\noindent 1) Assume $T$ is imaginary. If $p$ is not foreign to $\PP$, then for $a\models p$ there
exists $a'\in dcl(Aa)\backslash acl(A)$ such that $tp(a'/A)$ is $\PP$-internal.

\noindent 2) Assume $T$ is imaginary. Then $p$ is analyzable in $\PP$ iff every non-algebraic
extension of $p$ is non-foreign to $\PP$.

\noindent 3) For a general simple theory 1) and 2) are true in the hyperimaginary context (where
"non-algebraic" is replaced by "unbounded").
\end{fact}

An easy fact we will be using is the following.

\begin{fact}\label{some internal}
Assume $tp(a_i)$ are $\PP$-internal for $i<\alpha$. Then $tp(\langle a_i \vert i<\alpha \rangle)$
is $\PP$-internal.
\end{fact}

An important characterization of almost-internality is the following fact [S0, Theorem 5.6.] (a
similar result obtained independently in [W, Proposition 3.4.9]).

\begin{fact}\label{finite-cover}
Let $p \in S(A)$ be an amalgamation base and let $\UU$ be an $A$-invariant set. Suppose $p$ is
almost-$\UU$-internal. Then there is a Morley sequence $\bar a$ in $p$ and there is a definable
relation $R(x,\bar y,\bar a)$ (over $\bar a$ only) such that, for every tuple $\bar c$, $R(\CC,\bar
c,\bar a)$ is finite and for every $a'$ realizing $p$, there is some tuple $\bar c$ from $\UU$ such
that $R(a',\bar c,\bar a)$ holds.
\end{fact}

\noindent $T$ is said to be \em unidimensional \em if whenever $p$ and $q$ are complete
non-algebraic types, $p$ and $q$ are non-orthogonal. An $A$-invariant set $\UU$ is called
supersimple if $SU(a/A)<\infty$ for every $a\in\UU$. From Fact \ref{finite-cover} and Fact
\ref{internal} it is easy to deduce the following (using compactness).

\begin{fact}\label{supersimple definable}
Let $T$ be a simple theory. Let $p\in S(\emptyset)$ and let $\theta\in L$. Assume $p$ is analyzable
in $\theta^\CC$. Then $p$ is analyzable in $\theta^\CC$ in finitely many steps. In particular, if
$T$ is an imaginary simple unidimensional theory and there exists a non-algebraic supersimple
definable set, then $T$ has finite $SU$-rank, i.e. every complete type has finite $SU$-rank (in
fact, for every given sort there is a finite bound on the $SU$-rank of all types in that sort).
\end{fact}

\subsection {The forking topology}


\begin{definition} \label{tau definition}
\em Let $A\subseteq \CC$. An invariant set $\UU$ over $A$ is said to be \em a basic $\tau^f$-open
set over $A$ \em if there is $\phi(x,y)\in L(A)$ such that $$\UU=\{a \vert \phi(a,y)\ \mbox{forks\
over}\ A \}.$$
\end{definition}

\noindent Note that the family of basic $\tau^f$-open sets over $A$ is closed under finite
intersections, thus form a basis for a unique topology on $S_x(A)$.

\begin{definition}\label{projection closed}
\em We say that \em the $\tau^f$-topologies over $A$ are closed under projections ($T$ is PCFT over
$A$) \em if for every $\tau^f$-open set $\UU(x,y)$ over $A$ the set $\exists y \UU(x,y)$ is a
$\tau^f$-open set over $A$. We say that \em the $\tau^f$-topologies are closed under projections
($T$ is PCFT) \em if they are over every set $A$.
\end{definition}

We will make an essential use of the following facts from [S2].

\begin{fact}\label{tau extensions}
Let $\UU$ be a $\tau^f$-open set over a set $A$ and let $B\supseteq A$ be any set. Then $\UU$ is
$\tau^f$-open over $B$.
\end{fact}

We say that an $A$-invariant set $\UU$ has $SU$-rank $\alpha$ and write $SU(\UU)=\alpha$ if
$Max\{SU(p) \vert p\in S(A), p^\CC\subseteq\UU\}=\alpha$. We say that an $A$-invariant set $\UU$
has bounded finite $SU$-rank if there exists $n<\omega$ such that $SU(\UU)=n$.

\begin{fact}\label{tau bounded SU}
Let $\UU$ be an unbounded $\tau^f$-open set over some set $A$. Assume $\UU$ has bounded finite
$SU$-rank. Then there exists a set $B\supseteq A$ and $\theta(x)\in L(B)$ of $SU$-rank 1 such that
$\theta^\CC\subseteq \UU\cup acl(B)$.
\end{fact}

The following theorem [S2, Theorem 3.11] generalizes Fact \ref{supersimple definable} but at the
price of PCFT.

\begin{fact}\label{main tau}
Assume $T$ is a simple theory with PCFT. Let $p\in S(A)$ and let $\UU$ be a $\tau^f$-open set over
$A$. Suppose $p$ is analyzable in $\UU$. Then $p$ is analyzable in $\UU$ in finitely many steps.
\end{fact}

\section {Unidimensionality and PCFT}

In [BPV] it is defined when in a simple theory the extension property is first-order. Pillay [P1]
proved, using the result on the elimination of the $\exists^\infty$ [S1], that in any
unidimensional simple theory the extension property is first-order. Here we show that if $T$ is any
simple theory in which the extension property is first-order then $T$ is PCFT. We conclude that any
unidimensional simple theory is PCFT. From this we obtain, by Fact \ref{main tau}, the first step
towards the main result, namely, the existence of an unbounded $\tau^f$-open set that is
supersimple in an imaginary simple unidimensional theory implies $T$ is supersimple. In this
section $T$ is assumed to be simple, and if not stated otherwise, we work in $\CC$, however we
start with some notions that we will need for
hyperimaginaries.\\

First, we introduce some natural extensions of notions from [BPV]. By a pair $(M,P^M)$ of $T$ we
mean an $L_P=L\cup \{P\}$-structure, where $M$ is a model of $T$ and $P$ is a new predicate symbol
whose interpretation is an elementary submodel of $M$.  For the rest of this section, by a $\vert T
\vert$-\em small type \em we mean a complete hyperimaginary type in $\leq \vert T\vert$ variables
over a hyperimaginary of length $\leq \vert T\vert$.

\begin{definition}\label {def1}
Let $\PP_0,\PP_1$ be $\emptyset$-invariant families of $\vert T\vert$-small types.

\noindent 1) We say that a pair $(M,P^M)$ satisfies the extension property for $\PP_0$ if for every
$L$-type $p\in S(A),\ A\in dcl(M)$ with $p\in \PP_0$ there is $a\in p^M$ such that
$\nonfork{a}{P^M}{A}$.

\noindent 2) Let $$T_{Ext,\PP_0}=\bigcap \{Th_{L_P}(M,P^M) \vert\ \mbox{the\ pair\ $(M,P^M)$
satisfies\ the\ extension\ property\  w.r.t.\  } \PP_0\ \}.$$

\noindent 3) We say that $\PP_0$ dominates $\PP_1$ w.r.t. the extension property if $(M,P^M)$
satisfies the extension property for $\PP_1$ for every $\vert T\vert^+$-saturated pair
$(M,P^M)\models T_{Ext,\PP_0}$. In this case we write $\PP_0\unrhd_{_{Ext}} \PP_1$.

\noindent 4) We say that the extension property is first-order for $\PP_0$ if
$\PP_0\unrhd_{_{Ext}}\PP_0$. We say that the extension property is first-order if the extension
property is first-order for the family of all $\vert T\vert$-small types (equivalently, for the
family of all real types over sets of size $\leq \vert T\vert$).
\end{definition}

\begin {fact} $[BPV, Proposition\ 4.5]$\label {foext}
The extension property is first-order in $T$ iff for every formulas $\phi(x,y),\psi(y,z)\in L$ the
relation $Q_{\phi,\psi}$ defined by: $$Q_{\phi,\psi}(a)\mbox{\ iff}\ \phi(x,b)\mbox{ doesn't\ fork\
over}\ a\ \mbox{for\ every}\ b\models \psi(y,a)$$ is type-definable (here $a$ can be an infinite
tuple from $\CC$ whose sorts are fixed).
\end {fact}

Now, recall the following two facts and their corollary. First, let $\PP^{SU\leq 1}$ denote the
class of complete types over sets of size $\leq \vert T \vert$, of $SU$-rank $\leq 1$.

\begin{fact} \label {fact 1} $[P1]$
Let $T$ be a simple theory that eliminates $\exists^\infty$. Moreover, assume every non-algebraic
type is non-foreign to $\PP^{SU\leq 1}$. Then the extension property is first-order in $T$.
\end{fact}

\begin{fact} $[S1]$
Let $T$ be any unidimensional simple theory. Then $T$ eliminates $\exists^\infty$.
\end{fact}

\begin{corollary}\label{uni ext prop}
In any unidimensional simple theory the extension property is first-order.
\end{corollary}

Here we give an easy generalization of Fact \ref{fact 1}. For an $\emptyset$-invariant family
$\PP_0$ of $\vert T\vert$-small types we say that $\PP_0$ \em is extension-closed \em if for all
$p\in \PP_0$ if $\bar p$ is any extension of $p$ to a $\vert T\vert$-small type, then $\bar p\in
\PP_0$. First, we need an easy remark.

\begin{remark}\label{dcl T}
1) Assume $a$ is a hyperimaginary of length $\leq \vert T\vert$, and $B$ is a small set of
hyperimaginaries. Assume $a\in dcl(B)$. Then there exists $B_0\subseteq B$ of size $\leq \vert
T\vert$ such that $a\in dcl(B_0)$.\\
2) If $a$ is a hyperimaginary of length $\leq\vert T\vert$ and $b\in dcl(a)$ is arbitrary
hyperimaginary then $b$ interdefinable with a hyperimaginary of length $\leq\vert T\vert$.

\end{remark}

\proof 1) First, note there are hyperimaginaries of countable length $a_i=\bar a_i/E_i$, for $i\in
\vert T\vert$, where each $E_i$ is a type-definable equivalence relation over $\emptyset$ that
consists of countably many formulas such that $a$ is interdefinable with $(a_i\vert i\in \vert
T\vert)$ (by repeated applications of compactness). Thus we may assume that $a\in dcl(B)$ and
$a=\bar a/E$ where the length of $\bar a$ is countable and $E$ consists of countably many formulas.
Indeed, assuming this, we get that $tp(a/B)\vdash E(x,\bar a)$. Thus, by compactness, for every
formula $\psi(x)\in E(x,\bar a)$ there is a formula $\phi(x)\in tp(a/B)$, such that $\phi(x)\vdash
\psi(x)$, in particular there is a countable $B_0\subseteq B$ such that $tp(a/B_0)\vdash E(x,\bar
a)$. Hence $a\in dcl(B_0)$. 2) is easy and left to the reader.

\begin{lemma} \label{ext_family}
Let $\PP_0$ be an $\emptyset$-invariant family of $\vert T\vert$-small types. Assume $\PP_0$ is
extension-closed and that the extension property is first-order for $\PP_0$. Let $\PP^*$ be the
maximal class of $\vert T\vert$-small types such that $\PP_0\unrhd_{_{Ext}}\PP^*$. Then
$\PP^*\supseteq An(\PP_0)$, where $An(\PP_0)$ denotes the class of all $\vert T\vert$-small types
analyzable in $\PP_0$ by hyperimaginaries.
\end{lemma}

\proof Note that if the pair $(M,P^M)$ satisfies the extension property for the family of
$\emptyset$-conjugates of a hyperimaginary type $tp(b/A)$ and for the family of
$\emptyset$-conjugates of some hyperimaginary type $tp(a/bA)$ then $(M,P^M)$ satisfies the
extension property for the family of $\emptyset$-conjugates of $tp(ab/A)$. Thus, since $\PP_0$ is
extension-closed and the extension property is first-order for $\PP_0$ we conclude that if $B$ is
any hyperimaginary of length $\leq \vert T\vert$ and $\bar a$ is a tuple of of length $\leq \vert
T\vert$ of realizations of some types from $\PP_0$ over $B$, then if $(M,P^M)$ is a $\vert T
\vert^+$-saturated pair and $(M,P^M)\models T_{Ext,\PP_0}$ then $(M,P^M)$ satisfies the extension
property for the family of $\emptyset$-conjugates of $tp(\bar a/B)$. Now, assume $tp(a/A)$ is a
$\vert T\vert$-small type that is $\PP_0$-internal. There is a set $B$ with $A\in dcl(B)$ such that
$a$ is independent from $B$ over $A$ and there is a tuple of realizations $\bar c$ of types from
$\PP_0$ over $B$ such that $a\in dcl(B\bar c)$. By Remark \ref{dcl T}(1), we may assume both $B$
and $\bar c$ are of length $\leq\vert T\vert$. By the previous observation, $tp(\bar c/B)\in
\PP^*$. Since $a\in dcl(B\bar c)$, and $a$ is independent from $B$ over $A$ we conclude
$tp(a/A)\in\PP^*$. Now, assume $tp(a/A)$ is a $\vert T\vert$-small type that is analyzable in
$\PP_0$ by hyperimaginaries. By repeated applications of Fact \ref{internal} in the hyperimaginary
context and Remark \ref{dcl T}(2), for some $\alpha<\vert T\vert^+$ there exists a sequence $(a_i
\vert i\leq\alpha)\subseteq dcl(aA)$ of hyperimaginaries of length $\leq\vert T\vert$ such that
$a_\alpha=a$ and such that $tp(a_i/\{a_j \vert j<i\}\cup A)$ is $\PP_0$-internal for every
$i\leq\alpha$. By the previous observation $tp(a_i/\{a_j \vert j<i\}\cup A)\in \PP^*$ for every
$i\leq\alpha$. By applying the first observation inductively we get that $tp((a_i \vert
i\leq\alpha)/A)\in \PP^*$, and in particular $tp(a/A)\in \PP^*$.

\begin{remark}
Note that if $T$ eliminates $\exists^\infty$ then the extension property is first-order for
$\PP^{SU\leq 1}$ (this was proved in [V, Proposition 2.15]). Thus Lemma \ref{ext_family} implies
Fact \ref{fact 1}.
\end{remark}

Now, we aim to show that any simple theory in which the extension property is first-order is PCFT.

\begin{definition}
We say that $T$ is semi-PCFT over $A$ if for every formula $\psi(x,yz)\in L(A)$ the set $\{a \vert
\ \psi(x,ab) \mbox{\ forks\ over\ } Aa\  \mbox{for\ some\ } b \}$ is $\tau^f$-open over $A$.
\end{definition}

\begin{lemma}\label{ext over sets}
\noindent 1) If the extension property is first-order then the extension property is first-order
over every set $A$.\\
\noindent 2) If the extension property is first-order, then $T$ is semi-PCFT over $\emptyset$.
Thus, if the extension property is first-order, then $T$ is semi-PCFT over every set $A$.
\end{lemma}

\noindent

\proof 1) Let $\phi(x,y,A),\psi(y,z,A)\in L(A)$. Let $Q'$ be the relation defined by $Q'(a'A')$ iff
$\phi(x,bA')$ doesn't fork over $a'A'$ for all $b\models \psi(y,a'A')$. Clearly, for all $a'A'$ we
have $Q'(a'A')$ iff $\phi(x,bA'')$ doesn't fork over $a'A'$ for all $bA''$ such that $b\models
\psi(y,a'A')$ and $A''=A'$ (of course, $A''$ can be taken to be a finite tuple and we only need to
require that $A''$ is equal to certain coordinates of $A'$). By Fact \ref{foext} we see that $Q'$
is type-definable. In particular, $\{ a' \vert\ Q(a'A)\}$ is type-definable over $A$. Thus by Fact
\ref{foext}, the extension property if first-order over $A$.

\noindent 2) Assume the extension property is first-order. Let $\psi(x,yz)\in L$, we need to show
that the set $F=\{a \vert \ \psi(x,ab) \mbox{\ doesn't\ fork\ over\ } a\ \mbox{for\ all\ } b \}$ is
$\tau^f$-closed. Indeed, clearly $$F=\{a \vert \ \psi(x,a'b) \mbox{\ doesn't\ fork\ over\ } a\
\mbox{for\ all\ } a'b\ \mbox{with\ }a'=a \}.$$ By Fact \ref{foext}, $F$ is Stone-closed, in
particular $F$ is $\tau^f$-closed.

\begin{lemma}
Assume $T$ is semi-PCFT over $A$. Then $T$ is PCFT over $A$.
\end{lemma}

\proof We may clearly assume $A=\emptyset$. Let $\psi(x,yz)\in L$. We need to show that $\Gamma$,
defined by $\Gamma(a)$ iff $\forall b (\psi (x,ab) \mbox{\ doesn't\ fork\ over}\ \emptyset)$ is a
$\tau^f$-closed set. Let $\Gamma^*$ be defined by: for all $a$: $$\Gamma^*(a)\mbox{\ \ iff\ }
\bigwedge_{\phi(x,y)\in L} [\phi(x,a)\mbox{\ forks\ over\ }\emptyset\rightarrow \forall b
(\psi(x,ab)\wedge \neg\phi(x,a) \mbox{\ doesn't\ fork\ over\ }a)].$$

To finish it is sufficient to prove:

\begin{subclaim}
$\Gamma^*$ is $\tau^f$-closed and $\Gamma=\Gamma^*$.
\end{subclaim}

\proof First, by our assumption $\Gamma^*$ is $\tau^f$-closed. To prove the second part, first
assume $\Gamma(a)$. Then for any $b$ there is $c$ such that $\nonforkempty{c}{ab}$ and
$\psi(c,ab)$. Thus $\Gamma^*(a)$. Assume now $\Gamma^*(a)$. Let $p^{ind}_a(x)=\bigwedge { \{
\neg\phi(x,a) \vert\ \phi(x,y)\in L,\  \phi(x,a)\mbox{\ forks\ over\ }\emptyset \} }$. Let $b$ be
arbitrary and let $q(x)=p^{ind}_a(x)\wedge \psi(x,ab)$. It is enough to show that $q(x)$ doesn't
fork over $a$ (since any realization of $q$ is independent of $a$). Indeed, by $\Gamma^*(a)$, every
finite subset of $q(x)$ doesn't fork over $a$, so we are done.

\begin{corollary}\label{ext pcft}
Suppose the extension property is first-order in $T$. Then $T$ is PCFT.
\end{corollary}

Combining the last two corollaries we get:

\begin{theorem}\label {uni_pcft}
Let $T$ be any unidimenisonal simple theory. Then $T$ is PCFT.
\end{theorem}

\begin{corollary}\label {ss_tau_sets}
Let $T$ be an imaginary simple unidimensional theory. Let $p\in S(A)$ and let $\UU$ be an unbounded
$\tau^f$-open set over $A$. Then $p$ is analyzable in $\UU$ in finitely many steps. In particular,
for such $T$ the existence of an unbounded supersimple $\tau^f$-open set over some set $A$ implies
$T$ is supersimple.
\end{corollary}

\proof By Theorem \ref{uni_pcft} every unidimensional theory is $PCFT$. Thus by Fact \ref{main tau}
and the assumption that $T$ is imaginary and unidimensional, if $\UU$ is an unbounded $\tau^f$-open
set over $A$, then $tp(a/A)$  is analyzable in $\UU$ in finitely many steps for every $a\in\CC$.
Thus, if $\UU$ is supersimple, $SU(a/A)<\infty$ for all $a\in \CC$. Thus $T$ is supersimple.


\begin{remark}\label {tau comp}
Note that at this point we can conclude that any countable imaginary low unidimensional theory is
supersimple. Indeed, by Fact \ref{tau bounded SU} and Fact \ref{supersimple definable} it will be
sufficient to show the existence of an unbounded $\tau^f$-open set of bounded finite $SU$-rank. The
existence of such a set follows by Hrushovski's Baire categoricity argument [H0] together with
Theorem \ref{uni_pcft} applied to the $\tau^f$-topology (see [P]). The reason this argument works
is that the $\tau^f$-topology in a low theory is a Baire space. Indeed, in a low theory basic
$\tau^f$-open sets are type-definable and therefore we have the following property: the
intersection of a countable chain of basic $\tau^f$-open sets is non-empty iff each set in the
chain is non-empty. If the theory is not low we don't know the above property is true for the basic
$\tau^f$-open sets.
\end{remark}

\section {Definability of being in the canonical base}
In this section we show that in suitable setting the relation $R$ defined by $R(e,a)$ iff $e\in
acl(Cb(C/a))$ is Stone open over $C$ for a fixed set $C$. This definability result will be crucial
for the dichotomy theorem we will prove in the next section. In this section $T$ is assumed to be a
simple theory and all tuples and sets are assumed to be from $\CC$, however $Cb(A/B)$, for sets
$A,B$, is the canonical base of $Lstp(A/B)$ given as a hyperimaginary.

\begin{definition}
Let $C$ be any set. We say that a set $\UU$ is a basic $\tau^f_*$-open set over $C$ if there exists
$\psi(x,y,C)\in L(C)$ such that $\UU=\{a \vert \ \psi(x,aC) \mbox{\ forks\ over } a\}$.
\end{definition}

First, we note the following claim:

\begin{claim}\label {lemma tau star}
For every $e,C,a$, we have $e\in acl(Cb(C/a))$ iff for every Morley sequence $(C_i\vert i<\omega)$
of $Lstp(C/a)$ we have $e\in acl(C_i \vert i<\omega)$.
\end{claim}

\proof Left to right follows from the well known fact that $Cb(C/a)\in dcl(C_i \vert i<\omega)$ for
every Morley sequence $(C_i\vert i<\omega)$ of $Lstp(C/a)$. For the other direction, assume the
right hand side. Let $(C_i \vert i<\omega\cdot 2)$ be a Morley sequence of $Lstp(C/a)$. Let
$e^*=Cb(C/a)$. Then $e^*\in bdd(a)$ and thus clearly $(C_i \vert i<\omega\cdot 2)$ is a Morley
sequence of $Lstp(C/e^*)$. In particular, $(C_i \vert i<\omega)$ is independent from $(C_i \vert
\omega\leq i<\omega\cdot 2)$ over $e^*$. By our assumption, $e\in acl(C_i \vert i<\omega)$ and
$e\in acl(C_i \vert \omega\leq i<\omega\cdot 2)$. Thus $e\in acl(e^*)$.

\begin{lemma}\label {tau star}
Let $C$ be any set and let $\WW=\{(e,a) \vert\ e\in acl(Cb(C/a))\}$ (where $e,a$ are taken from
fixed sorts). Then $\WW$ is a $\tau^f_*$-open set over $C$.
\end{lemma}

\proof First note that since $T$ is simple, for any two sorts, if $x,x'$ has the first sort, and
$y$ has the second sort, there exists a type-definable relation $E_L(x,x',y)$ such that for all
$a,a',b$ with the right sorts we have $E_L(a,a',b)$ iff $Lstp(a/b)=Lstp(a'/b)$. By Claim \ref{lemma
tau star}, $(e,a)\not\in \WW$ iff there exists an $a$-indiscernible sequence $(C_i\vert i<\omega)$
which is independent over $a$ with $E_L(C_0,C,a)$ such that $e\not\in acl(C_i \vert i<\omega)$. For
each $n<\omega$, let $$L_n=\{\bar\psi \vert\ \bar\psi=\{\psi_i(Y_0,Y_1,...,Y_i,y) \vert i\leq
n\}\mbox{\ for\ some\ } \psi_i\in L\}$$ (where $y$ has the sort of the $a$-s in $\WW$, and each
$Y_i$ has the sort of $C$). For each $n<\omega$ and $\bar\psi=\{\psi_i(Y_0,Y_1,...,Y_i,y) \vert
i\leq n\}\in L_n$, let
$\Theta_{\bar \psi}(x,y,Y_0,Y_1,...,Y_n,C)=$ \\
$$E_L(Y_0,C,y)\wedge I(Y_0,...Y_n,y)\wedge(\bigwedge_{i=0}^n
\neg\psi_i(Y_0,Y_1,...,Y_i,y))\wedge x\not\in acl(Y_0,Y_1,...,Y_n),$$ where $I(Y_0,...Y_n,y)$ is
the partial type saying $Y_0,...Y_n$ is $y$-indiscernible. Note that each $\Theta_{\bar
\psi}(x,y,Y_0,Y_1,...,Y_n,C)$ is a type-definable relation over $C$. By compactness, $(e,a)\not\in
\WW$ iff $$\bigwedge_{\bar\psi=\{\psi_i\}_i\in L_n, n<\omega} [(\bigwedge_{i=0}^n
\psi_i(Y_0,Y_1,...Y_{i-1},C,a)\mbox{\ forks\ over\ } a)\rightarrow \exists Y_0,...Y_n\Theta_{\bar
\psi}(e,a,Y_0,Y_1,...,Y_n,C)].$$ We see that the complement of $\WW$ is an intersection of
$\tau_*^f$-closed sets over $C$  (clearly, every Stone-closed set over $C$ is $\tau_*^f$-closed
over $C$.)

\begin{proposition}\label {open Cb}
Let $q(x,y)\in S(\emptyset)$ and let $\chi(x,y,z)\in L$ be such that $\models \forall y\forall
z\exists^{<\infty} x\chi(x,y,z)$. Then the set $$\UU=\{(e,c,b,a) \vert\ e\in acl(Cb(cb/a))\}$$ is
relatively Stone-open inside the Stone-closed set
$$F=\{(e,c,b,a)\vert\ \nonforkempty{b}{a}, \models\chi(c,b,a), tp(cb)=q\}.$$ (where
$e$ is taken from a fixed sort too).
\end{proposition}

\proof Note that since $q\in S(\emptyset)$, it is enough to show that for any fixed $c^*b^*\models
q$ the set $\UU^*=\{(e,a) \vert\ e\in acl(Cb(c^*b^*/a))\}$ is relatively stone-open inside
$$F^*=\{(e,a)\vert\ \nonforkempty{b^*}{a}, \models\chi(c^*,b^*,a)\}.$$ Now, by Lemma \ref{tau star}
we know $\UU^*$ is a $\tau_*^f$-open set over $b^*c^*$. Thus, for some $\psi_i(t_i;w,z,c^*b^*)\in
L(c^*b^*)$ ($i\in I$) we have $\UU^*=\bigcup_i \UU^*_{\psi_i}$ where $\UU^*_{\psi_i}=\{(e,a)\vert\
\psi_i(t_i;e,a,c^*b^*) \mbox{\ forks over } ea\}$.

\begin{subclaim}\label{subclaim 1}
For every $(e,a)\in F^*$ we have $(e,a)\in \UU^*_{\psi_i}$ iff $$\forall d
(\psi_i(d;e,a,c^*b^*)\rightarrow \forkempty{da}{b^*})\wedge e\in acl(a).$$
\end{subclaim}

\proof Let $(e,a)\in F^*$. Assuming the left hand side we know $e\in acl(Cb(c^*b^*/a))$, hence
$e\in acl(a)$. Let $d\models\psi_i(z;e,a,c^*b^*)$. If $\nonforkempty{da}{b^*}$, then
$\nonfork{d}{b^*}{a}$. Since $(e,a)\in F^*$, $c^*\in acl(b^*a)$ implies $\nonfork{d}{b^*c^*}{ea}$,
contradicting $(e,a)\in \UU^*_{\psi_i}$. Assume now the right hand side. By a way of contradiction
assume there exists $d\models\psi_i(t_i;e,a,c^*b^*)$ such that $\nonfork{d}{b^*c^*}{ea}$. Since
$e\in acl(a)$, this equivalent to $\nonfork{d}{b^*c^*}{a}$. Since
$(e,a)\in F^*$ this is equivalent to $\nonforkempty{da}{b^*}$, contradiction.\\

\noindent By Subclaim \ref{subclaim 1} we see that each of $\UU^*_{\psi_i}$ and hence $\UU^*$ is
Stone-open relatively inside $F^*$.

\section {A dichotomy for projection closed topologies}
The main obstacle for proving that a countable imaginary simple unidimensional theory is
supersimple is, as indicated in Remark \ref{tau comp}, the lack of compactness. The goal of this
section is to prove a dichotomy that will enable us to reduce the general situation to a context
where compactness can be applied eventually. More specifically, we consider a general family of
topologies on the Stone spaces $S_x(A)$ that refine the Stone topologies and are closed under
projections (and under adding dummy variables). For any such family of topologies the dichotomy
says that either there exists an unbounded invariant set $\UU$ that is open in this topology and is
supersimple, OR for any $SU$-rank 1 type $p_0$ every type analyzable in $p_0$ is analyzable in
$p_0$ by essentially 1-based types by mean of our family of topologies. In this section $T$ is
assumed to be an imaginary simple theory and we work in $\CC=\CC^{eq}$.

\begin{definition}\em
A family $$\Upsilon=\{\Upsilon_{x,A} \vert\ x \mbox{ is a finite sequence of variables and }
A\subset \CC \mbox{ is small}\}$$ is said to be \em a projection closed family of topologies \em if
each $\Upsilon_{x,A}$ is a topology on $S_x(A)$ that refines the Stone-topology on $S_x(A)$, this
family is invariant under automorphisms of $\CC$ and change of variables by variables of the same
sort, and the family is closed under product by the full Stone space $S_y(A)$ (where $y$ is a
disjoint tuple of variables) and closed by projections, namely whenever $\UU(x,y)\in
\Upsilon_{xy,A}$, $\exists y\UU(x,y)\in\Upsilon_{x,A}$.
\end{definition}

There are two natural examples of projections-closed families of topologies; the Stone topology and
the $\tau^f$-topology of a PCFT theory. From now on fix a projection closed family $\Upsilon$ of
topologies.

\begin{definition}\label {def ess-1-based}\em
1) A type $p\in S(A)$ is said to be \em essentially 1-based over $A_0\subseteq A$, by mean of
$\Upsilon$ \em if for every finite tuple $\bar c$ from $p$ and for every type-definable
$\Upsilon$-open set $\UU$ over $A\bar c$, with the property that $a$ is independent from $A$ over
$A_0$ for every $a\in \UU$, the set $\{a\in \UU \vert\ Cb(a/A\bar c)\not\in bdd(aA_0)\}$ is nowhere
dense in the Stone-topology of $\UU$. We say $p\in S(A)$ is \em essentially 1-based by mean of
$\Upsilon$ \em if $p$ is essentially 1-based over $A$ by mean of $\Upsilon$.\\ 
2) Let $V$ be an $A_0$-invariant set and let $p\in S(A_0)$. We say that $p$ is \em analyzable in
$V$ by essentially 1-based types by mean of $\Upsilon$ \em if there exists $a\models p$ and there
exists a sequence $(a_i\vert\ i\leq\alpha)\subseteq dcl^{eq}(A_0a)$ with $a_\alpha=a$ such that
$tp(a_i/A_0\cup\{a_j\vert j<i\})$ is $V$-internal and essentially 1-based over $A_0$ by mean of
$\Upsilon$ for all $i\leq\alpha$.
\end{definition}


\begin{remark}
Note that $p\in S(A)$ is essentially 1-based by mean of $\Upsilon$ iff for every finite tuple $\bar
c$ from $p$ and for every non-empty type-definable $\Upsilon$-open set $\UU$ over $A\bar c$, there
exists a non-empty relatively Stone-open and Stone-dense subset $\chi$ of $\UU$ such that
$\nonfork{a}{\bar c}{acl^{eq}(Aa)\cap acl^{eq}(A\bar c)}$ for all $a\in \chi$. Intuitively,  $p\in
S(A)$ is essentially 1-based over a proper subset $A_0$ of $A$ by mean of $\Upsilon$ if the
canonical base of $Lstp(a/A\bar c)$ can be pushed down to $bdd(aA_0)$ for "most" $a\in\UU$ provided
that $a$ independent from $A$ over $A_0$ for all $a\in\UU$. This will be important for the
reduction in section 7.
\end{remark}

\begin{example}
The unique non-algebraic 1-type over $\emptyset$ in algebraically closed fields is not essentially
1-based by mean of the $\tau^f$-topologies.
\end{example}

\proof Work in a saturated algebraically closed field $\bar K$. Let $k_0\leq \bar K$ denote the
prime field, and $acl$ denote the algebraic closure in the home sort. First, recall that for every
finite tuples $\bar a,\bar b,\bar c\subseteq \bar K$ we have $\nonfork{\bar a}{\bar b}{\bar c}$ iff
$tr.deg (k_0(\bar a,\bar c)/k_0(\bar c))=tr.deg(k_0(\bar a,\bar b,\bar c)/k_0(\bar b,\bar c))$.
Now, it is known that if $c_0,c_1,a\in \bar K$ are algebraically independent over $k_0$, and
$b=c_1a+c_0$ then $acl(ab)\cap acl(c_0c_1)=acl(\emptyset)$, $tr.deg(k_0(a,b)/k_0)=2$, and
$tr.deg(k_0(a,b,c_0,c_1)/k_0(c_0,c_1))=1$, and therefore $\fork{ab}{c_0c_1}{acl^{eq}(ab)\cap
acl^{eq}(c_0c_1)}$ ($acl$ can be replaced by $acl^{eq}$ since $\bar K$ eliminates imaginaries). Let
$p\in S(\emptyset)$ be the unique non-algebraic 1-type. Let us fix two algebraically independent
realizations $c_0,c_1$ of $p$. Let $\UU$ be defined by:
$$\UU=\{(a,b)\in \bar K^2 \vert a\not\in acl(c_0,c_1)\mbox{ and\ } b=c_1a+c_0\}.$$ Note that $\UU$
is a type-definable $\tau^f$-open set over $c_0c_1$. The above observation shows $\UU$ fails to
satisfy the requirement in Definition \ref{def ess-1-based}(1). Thus $p$ is not essentially-1-based
by mean of the $\tau^f$-topologies.\\

One of the key ideas for proving the main result is the following theorem. We say that an
$A$-invariant set $\UU$ has finite $SU$-rank if $SU(a/A)<\omega$ for every $a\in\UU$.

\begin{theorem}\label{dichotomy thm}
Let $T$ be a countable simple theory that eliminates hyperimaginaries. Let $\Upsilon$ be a
projection-closed family of topologies. Let $p_0$ be a partial type over $\emptyset$ of $SU$-rank
1. Then, either there exists an unbounded finite-$SU$-rank $\Upsilon$-open set over some countable
set, or every type $p\in S(A)$, with $A$ countable, that is internal in $p_0$ is essentially
1-based over $\emptyset$ by mean of $\Upsilon$. In particular, either there exists an unbounded
finite $SU$-rank $\Upsilon$-open set, or whenever $A$ is countable, $p\in S(A)$ and every
non-algebraic extension of $p$ is non-foreign to $p_0$, $p$ is analyzable in $p_0$ by essentially
1-based types by mean of $\Upsilon$.
\end{theorem}

\proof $\Upsilon$ will be fixed and we'll freely omit the phrase "by mean of $\Upsilon$". To see
the "In particular" part, work over $A$ and assume that every $p'\in S(A')$, with $A'\supseteq A$
countable, that is internal in $p_0$, is essentially 1-based over $A$. Indeed, assume $p\in S(A)$
is such that every non-algebraic extension of $p$ is non-foreign to $p_0$. Then, for $a\models p$
there exists $a'\in dcl^{eq}(Aa)\backslash acl^{eq}(A)$ such that $tp(a'/A)$ is $p_0$-internal and
thus essentially 1-based over $A$ by our assumption. Since $L$ and $Aa$ are countable so is
$dcl^{eq}(Aa)$ and thus by repeating this process we get that $p$ is analyzable in $p_0$ by
essentially 1-based types. We prove now the main part. Assume there exist a countable $A$ and $p\in
S(A)$ that is internal in $p_0$ and $p$ is not essentially 1-based over $\emptyset$. 
By Fact \ref{some internal}, we may assume there exists $d\models p$, and $b$ that is independent
from $d$ over $A$, and a finite tuple $\bar c\subseteq p_0$ such that $d\in dcl(Ab\bar c)$, and
there exists a type-definable $\Upsilon$-open set $\UU$ over $Ad$ such that $a$ is independent from
$A$ for all $a\in \UU$ and $\{a\in \UU \vert Cb(a/Ad)\not\subseteq acl^{eq}(a)\}$ is not nowhere
dense in the Stone-topology of $\UU$. So, since $\Upsilon$ refines the Stone-topology, by
intersecting it with a definable set, we may assume that $\{a\in \UU \vert Cb(a/Ad)\not\subseteq
acl^{eq}(a)\}$ is dense in the Stone-topology of $\UU$. Now, for each disjoint partition $\bar
c=\bar c_0\bar c_1$ and formula $\chi(\bar x_1,\bar x_0,y,z)\in L(A)$ such that (*)\ $\forall \bar
x_0,y,z \exists^{<\infty} \bar x_1 \chi(\bar x_1,\bar x_0,y,z)$, let
$$F_{\chi,\bar c_0,\bar c_1}=\{ a\in \UU \vert\ \exists b',\bar c'_0,\bar c'_1\ \mbox{s.t.}\
tp(b'\bar c'_0\bar c'_1/Ad)=tp(b\bar c_0\bar c_1/Ad)\ \mbox{and}\ a\ \mbox{is\ independent\ from}$$
$$b'\bar c'_0\bar c'_1\ \mbox{over}\ Ad\ \mbox{and} \models\chi(\bar c'_1,\bar c'_0,b',a)\
\mbox{and}\ a\ \mbox{is\ independent\ from}\ Ab'\bar c'_0\ \mbox{over}\ \emptyset\}.$$ Let
$\PP_{\bar c}$ be the (finite) set of partitions of $\bar c$ into two subsets. Note that since $d$
is independent from $b$ over $A$, any $a\in\UU$ is independent from $Ab'$ whenever
$tp(b'/Ad)=tp(b/Ad)$ and $\nonfork{a}{b'}{Ad}$. Thus, since $p_0$ is a partial type over
$\emptyset$ of $SU$-rank $\leq 1$ we have $$\UU=\bigcup_{(\bar c_0,\bar c_1)\in\PP_{\bar c},\
\chi\models (*)} F_{\chi,\bar c_0,\bar c_1}.$$ Note that since we are fixing the type of $b'\bar
c'_0\bar c'_1$ over $Ad$, the sets $F_{\chi,\bar c_0,\bar c_1}$ are type-definable over $Ad$. Since
$L$ and $A$ are countable, by the Baire category theorem for the Stone-topology of the closed set
$\UU$, there exists $(\bar c^*_0\bar, c^*_1)\in\PP_{\bar c}$ and there is $\chi^*\models (*)$ such
that $F_{\chi^*,\bar c^*_0,\bar c^*_1}$ has non-empty interior in the Stone-topology of $\UU$.
Thus, we may assume that $\UU$ is a type-definable $\Upsilon$-open set over $Ad$ such that $\{a\in
\UU \vert Cb(a/Ad)\not\subseteq acl^{eq}(a)\}$ is dense in the Stone-topology of $\UU$ and for
every $a\in \UU$ there exists $b'\bar c'_0\bar c'_1\models tp(b\bar c^*_0\bar c^*_1/Ad)$ that is
independent from $a$ over $Ad$ and such that $\models\chi^*(\bar c'_1,\bar c'_0,b',a)$ and $a$ is
independent from$Ab'\bar c'_0$ over $\emptyset$. Let us now define a set $V$ over $Ad$ by
$$V=\{(\bar c'_0,\bar c'_1,b',a',e') \vert\ \mbox{if}\ tp(b'\bar c'_0\bar c'_1/Ad)=tp(b\bar
c^*_0\bar c^*_1/Ad)\ \mbox{and}\ a'\ \mbox{is\ independent\ from}$$  $$b'\bar c'_0\bar c'_1\
\mbox{over}\ Ad\ \mbox{and}\ a'\ \mbox{is\ independent\ from}\ Ab'\bar c'_0\ \mbox{over}\
\emptyset\ \mbox{and}\ \models\chi^*(\bar c'_1,\bar c'_0,b',a')$$ $$\mbox{then}\ e'\in
acl(Cb(Ab'\bar c'_0\bar c'_1/a'))\}.$$ Let $$V^*=\{e' \vert \exists a'\in\UU\ \forall b',\bar
c'_0,\bar c'_1\ V(\bar c'_0,\bar c'_1,b',a',e')\}.$$

\begin{subclaim}
$V^*$ is a $\Upsilon$-open set over $Ad$.
\end{subclaim}

\proof By Proposition \ref{open Cb}, we see that $V$ is a Stone-open set over $Ad$. Note that
Stone-open sets are closed under the $\forall$ quantifier (indeed, if $\UU(x,y)$ is Stone-open,
then the complement of $\forall y\UU(x,y)$ is Stone-closed by compactness). Therefore, since the
$\Upsilon$ topology refines the Stone-topology and closed under product by a full Stone-space and
closed under projections, we conclude that $V^*$ is a $\Upsilon$-open set.

\begin{subclaim}
For appropriate sort for $e'$, the set $V^*$ is unbounded and has finite $SU$-rank over $Ad$.
\end{subclaim}

\proof First, note the following.

\begin{remark}\label {dcl_cb remark}
Assume $d\in dcl(c)$. Then $Cb(d/a)\in dcl(Cb(c/a))$ for all $a$.
\end{remark}

\noindent Let $a^*\in\UU$ be such that $Cb(a^*/Ad)\not\subseteq acl^{eq}(a^*)$. Then
$Cb(Ad/a^*)\not\subseteq acl^{eq}(Ad)$. By Remark \ref{dcl_cb remark}, there exists $e^*\not\in
acl^{eq}(Ad)$ such that $e^*\in acl^{eq}(Cb(Ab'\bar c'_0\bar c'_1/a^*))$ for all $b'\bar c'_0\bar
c'_1\models tp(b\bar c^*_0\bar c^*_1/Ad)$. In particular, $e^*\in V^*$. Thus, if we fix the sort
for $e'$ in the definition of $V^*$ to be the sort of $e^*$, then $V^*$ is unbounded. Now, let
$e'\in V^*$. Then for some $a'\in\UU$, $\models V(\bar c'_0,\bar c'_1,b',a',e')$ for all $b',\bar
c'_0,\bar c'_1$. By what we saw above, there exists $b'\bar c'_0\bar c'_1\models tp(b\bar c^*_0\bar
c^*_1/Ad)$ that is independent from $a'$ over $Ad$ such that $\models\chi^*(\bar c'_1,\bar
c'_0,b',a')$ and $a'\ \mbox{is\ independent\ from}\ Ab'\bar c'_0\ \mbox{over}\ \emptyset$. Thus, by
the definition of $V^*$, $e'\in acl(Cb(Ab'\bar c'_0\bar c'_1/a'))$. Since $Ab'$ is independent from
$a'$ over $\emptyset$, $tp(e')$ is almost-$p_0$-internal, and thus
$SU(e')<\omega$. In particular, $SU(e'/Ad)<\omega$.\\

\noindent Thus $V^*$ is the required set.

\section{Stable dependence}
We introduce the relation stable dependence and show it is symmetric. In this section $T$ is
assumed to be a complete theory unless otherwise stated, and we work in $\CC^{eq}$.

\begin{definition}
Let $a\in \CC$, $A\subseteq B\subseteq \CC$. We say that $a$ is stably-independent from $B$ over
$A$ if for every stable $\phi(x,y)\in L$, if $\phi(x,b)$ is over $B$ (i.e. the canonical parameter
of $\phi(x,b)$ is in $dcl(B)$) and $a'\models \phi(x,b)$ for some $a'\in dcl(Aa)$, then $\phi(x,b)$
doesn't divide over $A$. In this case we denote it by $\snonfork{a}{B}{A}$.
\end{definition}

We will need some basic facts from local stability [HP]. From now on we fix a stable formula
$\phi(x,y)$. A formula $\psi\in L(\CC)$ is said to be a $\phi$-formula over $A$ if it is a finite
boolean combination of instances of $\phi$, that is equivalent to a formula with parameters from
$A$. A complete $\phi$-type over $A$ is a consistent complete set of $\phi$-formulas over $A$.
$S_\phi(A)$ denotes the set of complete $\phi$-types over $A$. Note that if $M$ is a model then
every $p\in S_\phi(M)$ is determined by the set $\{\psi\in p \vert\ \psi=\phi(x,a)\ \mbox{or}\
\psi=\neg\phi(x,a)\ \mbox{for}\ a\in M\}$ (in fact, it is easy to see that every $\phi$-formula
over $M$ is equivalent to a $\phi$-formula whose parameters are from $M$). Recall the following
well known facts.

\begin{fact} \label{LS fact}
Let $\phi(x,y)\in L$ be stable. Then\\
1) \em [HP, Lemma 5.4(i)]\em\ For any model $M$, every $p\in S_\phi(M)$ is definable.\\
2) \em [HP, Lemma 5.5]\em\ Let $A$ be any set, let $p\in S(A)$, and let $M\supseteq A$ be a model.
Then there exists $q\in
S_\phi(M)$ that is consistent with $p$ and is definable over $acl^{eq}(A)$.\\
3) \em[HP, Lemma 5.8]\em\ Let $A=acl(A)$. Let $p\in S_\phi(A)$. Then for every model $M\supseteq
A$, there is a unique $\bar p\in S_\phi(M)$ that extends $p$ and such that $\bar p$ is definable
over $A$ (i.e. its $\phi$-definition is over $A$). Moreover, there is a canonical formula over $A$
that is the
definition of any such $\bar p$ over any such model $M$.\\
4) \em[HP, Lemma 5.9]\em\ Assume $p,q\in S_\phi(acl(A))$ are such that $p\vert A=q\vert A$. Then
there exists $\sigma\in Aut(\CC/A)$ such that $\sigma (p)=q$.
\end{fact}

The following definition is standard.

\begin{definition}
Let $p\in S_\phi(B)$ and let $A\subseteq B$. We say that $p$ doesn't fork over $A$ in the sense of
local stability (=LS) if for some model $M$ containing $B$ and some $\bar p\in S_\phi(M)$ that
extends $p$, $\bar p$ is definable over $acl(A)$.
\end{definition}

\begin{claim}\label{stable claim}
Let $T$ be simple. Let $\phi(x,y)\in L$ be stable. Assume $\nonfork{a}{b}{A}$ and
$\nonfork{a'}{b}{A}$ and $Lstp(a/A)=Lstp(a'/A)$. Then $\phi(a,b)$ iff $\phi(a',b)$.
\end{claim}

\proof By definition, $Lstp(a/A)=Lstp(a'/A)$ iff there exist $a_0=a,...,a_n=a'$ such that for every
$i<n$ there is an infinite $A$-indiscernible sequence containing $(a_i,a_{i+1})$. By extension,
transitivity, and symmetry, we may assume $n=1$ and $b$ is independent from $aa'$ over $A$. Let
$(c_i \vert i\in Z\backslash\{0\})$ ($Z$ denotes the integers) be an $A$-indiscernible sequence
such that $c_{-1}=a$ and $c_1=a'$. Since $I=(c_{-i}c_i \vert i\in\omega\backslash\{0\})$ is
$A$-indiscernible and $b$ is independent from $aa'$ over $A$, we may assume $I$ is indiscernible
over $Ab$. We claim $\phi(a,b)$ iff $\phi(a',b)$. Indeed, otherwise we get
$\phi(c_i,b)\leftrightarrow\phi(c_j,b)$ iff
$i,j$ have the same sign; a contradiction to stability of $\phi(x,y)$.\\

The following lemma is easy but important.

\begin{lemma}\label{LS lemma}
Assume $T$ is a simple theory in which Lstp=stp over sets and let $\phi(x,y)\in L$ be stable. Then
for all $a$ and $A\subseteq B\subseteq\CC$, $tp_\phi(a/B)$ doesn't fork over $A$ in the sense of LS
iff $tp_\phi(a/B)$ doesn't fork over $A$.
\end{lemma}

\proof Assume $p_\phi=tp_\phi(a/B)$ doesn't fork over $A$ in the sense of LS. Extend it to a
complete $\phi$-type $\bar p_\phi$ over a sufficiently saturated and sufficiently
strongly-homogeneous model $\cal M$ that is definable over $acl(A)$. If $tp_\phi(a/B)$ divide over
$A$, there is an $acl(A)$-indiscernible sequence $(B_i \vert i<\omega)\subseteq \cal M$ such that
if $p^\phi_{B_i}$ are the corresponding $acl(A)$-conjugates of $p_\phi$, then $\bigwedge_i
p^\phi_{B_i}=\emptyset$. By the uniqueness of non-forking extensions (in the sense of LS) of
complete $\phi$-types over algebraically closed sets (and the fact that $\cal M$ is sufficiently
strongly-homogeneous) we conclude that $\bar p_\phi$ extends each $p^\phi_{B_i}$, a contradiction.
For the other direction, assume $p_\phi=tp_\phi(a/B)$ doesn't fork over $A$. Let ${\cal M}\supseteq
B$ be a sufficiently saturated and sufficiently strongly homogeneous model. Let $\bar p\in S({\cal
M})$ be an extension of $p_\phi$ that doesn't fork over $A$.
Let $\psi(y,c)\in L({\cal M)}$ be the definition of $\bar p\vert\phi$ (where $c$ is the canonical
parameter of $\psi$). We claim that $c\in acl(A)$. Indeed, otherwise let $\sigma\in Aut({\cal
M}/acl(A))$ be such that $\sigma c\neq c$. So, $\bar p, \sigma(\bar p)$ have different
$\phi$-definitions, a contradiction to Claim \ref{stable claim}.

\begin{corollary}\label{cor LS 1}
Let $T$ be a simple theory in which Lstp=stp over sets. Then for all $a,A\subseteq B\subseteq\CC$
we have $\snonfork{a}{B}{A}$ iff $tp_\phi(a'/B)$ doesn't fork over $A$ in the sense of LS for every
stable $\phi(x,y)\in L$ and every $a'\in dcl(aA)$.
\end{corollary}

Given $a,A\subseteq B\subseteq\CC$, we will say that $tp(a/B)$ \em doesn't fork over $A$ in the
sense of LS \em if the right hand side of Corollary \ref{cor LS 1} holds.

\begin{lemma}\label {s-sym-trans}
Let $T$ be a simple theory in which Lstp=stp over sets. Then

\noindent 1) stable independence is a symmetric relation, that is, for all $a,b,A$ we have
$\snonfork{a}{Ab}{A}$ iff $\snonfork{b}{Aa}{A}$.

\noindent 2) For all $a, A\subseteq B\subseteq C$, if $\snonfork{a}{B}{A}$ and
$\snonfork{a}{C}{B}$, then $\snonfork{a}{C}{A}$. In fact, in any theory the same is true in the
sense of LS.
\end{lemma}

\proof To prove 1), first note the following.

\begin{subclaim}\label {s-sym subclaim}
Let $\phi(x,y)\in L$ be stable and let $a,a'\in\CC$ and let $A\subseteq\CC$. Assume
$tp_\phi(a/A)=tp_\phi(a'/A)$. Then $\phi(a,y)$ forks over $A$ iff $\phi(a',y)$ forks over $A$.
\end{subclaim}

\proof Otherwise, there are $p,q\in S(\CC)$, both extends $tp_\phi(a/A)=tp_\phi(a'/A)$, and do not
fork over $A$ such that $p$ represent $\phi(x,y)$ (namely, for some $b\in M$, $\phi(x,b)\in p$) and
$q$ doesn't represent $\phi(x,y)$. By Fact \ref{LS fact} (4), $(p\vert\phi)\vert acl(A)$ and
$(q\vert\phi)\vert acl(A)$ are $A$-conjugate. Let $\sigma\in Aut(\CC/A)$ be such that
$\sigma((p\vert\phi)\vert acl(A))=(q\vert\phi)\vert acl(A)$. Now, both $\sigma(p\vert \phi)$ and
$q\vert \phi$ extend $(q\vert\phi)\vert acl(A)$ and doesn't fork over $acl(A)$, and therefore by
Lemma \ref{LS lemma}, both doesn't fork over $acl(A)$ in the sense of LS. By Fact \ref{LS fact}
(3), $\sigma(p\vert
\phi)=q\vert\phi$, which is a contradiction.\\

\noindent We prove symmetry. Assume $\snonfork{a}{Ab}{A}$. To show $\snonfork{b}{Aa}{A}$, let
$\phi(x,y)\in L$ be stable such that $\phi(b',a')$ for some $b'\in dcl(Ab)$ and some $a'\in
dcl(Aa)$. Let $\tilde\phi(y,x)=\phi(x,y)$. By the assumption, $tp_{\tilde\phi}(a'/Ab)$ doesn't fork
over $A$ (in the usual sense), so there exists $a''\models tp_{\tilde\phi}(a'/Ab)$ such that
$\nonfork{a''}{Ab}{A}$. Let $(a''_i \vert i<\omega)$ be a Morley sequence of $tp(a''/Ab)$. Now,
$b'\models \bigwedge_{i<\omega} \phi(x,a''_i)$. Thus $\phi(x,a'')$ doesn't fork over $A$. By
Subclaim \ref{s-sym subclaim}, $\phi(x,a')$ doesn't fork over $A$. 2) is immediate by Corollary
\ref{cor LS 1} and the fact that the relation of being a non-forking extension in the LS sense is a
transitive relation on complete $\phi$-types (where $\phi$ is a fixed stable formula).

\section{An unbounded $\tau^f_{\infty}$-open set of bounded finite $SU_s$-rank is sufficient}

In this section we apply the dichotomy theorem from section 5 in order to reduce the problem on
supersimplicity of countable imaginary simple unidimensional theories to the problem of finding a
$\tau^f_\infty$-open set of finite $SU_s$-rank (over a finite set). In this section $T$ is an
imaginary simple theory. We work in $\CC^{eq}$.

\begin{definition}
1) For $a\in \CC$ and $A\subseteq \CC$ the $SU_s$-rank is defined by induction on $\alpha$: if
$\alpha=\beta+1$, then $SU_s(a/A)\geq\alpha$ if there exists $B\supseteq A$ such that
$\sfork{a}{B}{A}$ and $SU_s(a/B)\geq\beta$. For limit $\alpha$, $SU_s(a/A)\geq\alpha$ if
$SU_s(a/A)\geq\beta$ for all $\beta<\alpha$.

\noindent 2) Let $\UU$ be an $A$-invariant set. We write $SU_s(\UU)=\alpha$ (the $SU_s$-rank of
$\UU$ is $\alpha$) if $Max\{SU_s(p) \vert p\in S(A), p^\CC\subseteq\UU\}=\alpha$. We say that $\UU$
has bounded finite $SU_s$-rank if for some $n<\omega$, $SU_s(\UU)=n$. Note that the $SU_s$-rank of
$\UU$ might, a priori, depend on the choice of the set $A$ over which $\UU$ is invariant.
\end{definition}

\begin{definition}\em
\em The $\tau^f_\infty$-topology \em on $S(A)$ is the topology whose basis is the family of
type-definable $\tau^f$-open sets over $A$.
\end{definition}

\begin{lemma}\label{stable forking}
For $a\in\CC$ and $A\subseteq B\subseteq \CC$, assume $tp(a/B)$ doesn't fork over $acl(aA)\cap
acl(B)$ and $\fork{a}{B}{A}$. Then $\sfork{a}{B}{A}$.
\end{lemma}

\proof It will be sufficient to show that whenever $\fork{a}{B}{A}$  and $\nonfork{a}{B}{acl(a)\cap
acl(B)}$ for some (possibly infinite) tuple $a$ and some $A\subseteq B$, there exists a stable
$\phi(x,y)\in L$ such that $\phi(a,B)$ and $\phi(x,b)$ forks over $A$ (indeed, the above implies
the following: if $\fork{aA}{B}{A}$ and $\nonfork{aA}{B}{acl(aA)\cap acl(B)}$ then there exists a
stable formula $\phi(x,y)\in L$ such that $\phi(aA,B)$ and $\phi(x,B)$ forks over $A$, i.e.
$\sfork{a}{B}{A}$). To prove this, let $E=Cb(a/B)$. Then $E\subseteq acl(a)\cap acl(B)$. By the
assumption, there is $e^*\in dcl(E)\backslash acl(A)$, so $e^*\in (acl(a)\cap acl(B))\backslash
acl(A)$. Hence there are $n_0,n_1\in\omega$ and formulas $\chi_0(x,y),\chi_1(x,z)\in L$ such that
$\forall y\exists^{<n_0} x\chi_0(x,y)$ and $\forall z\exists^{<n_1} x\chi_0(x,z)$ and
$\chi_0(e^*,a)$ and $\chi_1(x,B)$ isolates $tp(e^*/B)$. Let $$\phi(y,z)\equiv \exists x
(\chi_0(x,y)\wedge \chi_1(x,z)).$$ Note that $\phi(y,z)$ is stable. Indeed, otherwise there are
$a\in \CC$ and an $\emptyset$-indiscernible sequence $B=(b_i \vert i\in Z)$ ($Z$=the integer
numbers) such that $i\geq 0$ iff $\phi(a,b_i)$. Since $B$ is indiscernible, and $\chi_1(x,b_0)$ is
algebraic, $\bigcap_{i\in I} \chi_1(\CC,b_i)=\bigcap_{i\in Z} \chi_1(\CC,b_i)$ for every infinite
$I\subseteq Z$. But since $\chi_0(x,a)$ is algebraic, for some infinite $I^*\subseteq \omega$,
$\chi_0(\CC,a)\cap \bigcap_{i\in I^*} \chi_1(\CC,b_i)\not=\emptyset$. A contradiction to
$\neg\phi(a,b_i)$ for $i<0$. To see that $\phi(y,B)$ forks over $A$, note that otherwise there
exists $a'\models \phi(y,B)$ such that $\nonfork{a'}{B}{A}$, so if $e'\models \chi_0(x,a')\wedge
\chi_1(x,B)$, then on one hand $\nonfork{e'}{B}{A}$ and on the other hand since $\chi_1(x,B)$
isolates $tp(e^*/B)$, $e'\in acl(B)\backslash acl(A)$ which is a contradiction.

\begin{lemma}\label {su_s finite to 1}
Assume $\UU$ is an unbounded $\tau^f_\infty$-open set of bounded finite $SU_s$-rank over some
finite set $A$. Then there exists a $\tau^f_\infty$-open set $\UU^*\subseteq \UU$ over some finite
set $B^*\supseteq A$ of $SU_s$-rank 1.
\end{lemma}

\proof We may clearly assume $\UU$ is a basic $\tau^f_\infty$-open set. Let $n=SU_s(\UU)$ ($\UU$ is
over $A$, and $n<\omega$). Let $a^*\in \UU$ with $SU_s(a^*/A)=n$. Let $B\supseteq A$ be finite such
that $\sfork{a^*}{B}{A}$, and $SU_s(a^*/B)=n-1$. So, there exists $a'\in dcl(a^*A)$ and stable
$\phi(x,y)\in L$ such that $\phi(a',B)$ and $\phi(x,B)$ forks over $A$. Let $f$ an
$\emptyset$-definable function such that $a'=f(a^*,A)$. Let $$\UU'=\{a\in\UU \vert\ \phi(f(a,A),B)\
\}\ (\mbox{as\ a\ set\ over\ }B).$$ Since $a^*\in \UU'$, $SU_s(\UU')\geq n-1$. If $a\in \UU'$, then
$\phi(f(a,A),B)$ implies $\sfork{a}{B}{A}$ and therefore $SU_s(\UU')\leq n-1$. We conclude
$SU_s(\UU')=n-1$. Clearly, $\UU'\subseteq\UU$ and $\UU'$ is type-definable. By Fact \ref{tau
extensions}, $\UU'$, is a $\tau^f$-open set over $B$. We finish by induction.

\begin{lemma}\label{main lemma}
Let $T$ be a countable imaginary simple unidimensional theory. Assume there is $p_0\in
S(\emptyset)$ of $SU$-rank 1 and there exists an unbounded $\tau^f_\infty$-open set over some
finite set of bounded finite $SU_s$-rank. Then $T$ is supersimple.
\end{lemma}

\proof  By Lemma \ref{su_s finite to 1}, there exists a finite set $A_0$ and a $\tau^f_\infty$-open
set $\UU$ over $A_0$ of $SU_s$-rank 1. Clearly, we may assume $\UU$ is type-definable. By Theorem
\ref{uni_pcft}, $T$ is PCFT. Thus, working over $A_0$, by Theorem \ref{dichotomy thm} for the
$\tau^f$-topology either (i) there exists an unbounded $\tau^f$-open set of finite $SU$-rank over
some countable set or (ii) every non-algebraic type over $A_0$ is analyzable in $p_0$ by
essentially 1-based types by mean of $\tau^f$. By Corollary \ref{ss_tau_sets}, we may assume (ii).
We claim $SU(\UU)=1$. Indeed, otherwise there exists $a$ and $d\in \UU$ such that
$\fork{d}{a}{A_0}$ and $d\not\in acl(aA_0)$. By (ii), there exists $(a_i \vert
i\leq\alpha)\subseteq dcl^{eq}(aA_0)$ with $a_\alpha=a$ such that $tp(a_i/A_0\cup\{a_j \vert
j<i\})$ is essentially 1-based over $A_0$ by mean of $\tau^f$ for all $i\leq\alpha$. Now, let
$i^*\leq\alpha$ be minimal such that there exists $d'\in\UU$ satisfying $\fork{d'}{\{a_i \vert
i\leq i^*\}}{A_0}$, and $d'\not\in acl(A_0\cup \{a_i \vert i\leq i^*\})$. Pick $\phi(x,a')\in
L(A_0\cup \{a_i \vert i\leq i^*\})$ that forks over $A_0$ and such that $\phi(d',a')$. Let
$$V=\{d\in\UU \vert\ \phi(d,a')\ \mbox{and}\ d\not\in acl(A_0\cup \{a_i \vert i\leq i^*\})\ \}.$$
By minimality of $i^*$, $d$ is independent from $\{a_i \vert i<i^*\}$ over $A_0$ for all $d\in V$.
Clearly $V$ is type-definable and by Fact \ref{tau extensions}, $V$ is a $\tau^f$-open set over
$A_0\cup \{a_i \vert i\leq i^*\}$. Now, since $tp(a_{i^*}/A_0\cup\{a_i \vert i<i^*\})$ is
essentially 1-based over $A_0$ by mean of $\tau^f$, the set
$$\{d\in V \vert\ Cb(d/A_0\cup \{a_i \vert i\leq i^*\})\in bdd(dA_0)\}$$ contains a relatively
Stone-open and Stone-dense subset of $V$. In particular, there exists $d^*\in V$ such that
$tp({d^*}/{A_0\cup \{a_i \vert i\leq i^*\}}$ doesn't fork over ${acl(A_0d^*)\cap acl(A_0\cup\{a_i
\vert i\leq i^*\})}$. Since we know $\fork{d^*}{A_0\cup\{a_i \vert i\leq i^*\}}{A_0}$, Lemma
\ref{stable forking} implies $\sfork{d^*}{A_0\cup\{a_i \vert i\leq i^*\}}{A_0}$. Hence $d^*\in V$
implies $SU_s(d^*/A_0)\geq 2$, which contradict $SU_s(\UU)=1$. Thus we have proved $SU(\UU)=1$.
Now, by Fact \ref{tau bounded SU} there exists a definable set of $SU$-rank 1, and thus by Fact
\ref{supersimple definable}, $T$ is supersimple.

\begin{remark}\label{baire remark}
Note that if $X$ is any Stone-closed subset of the Stone-space $S_x(T)$ and $B=\{F_i\}_{i\in I}$ is
a basis for a topology $\tau$ on $X$ that consists of Stone-closed subsets of $X$, then $(X,\tau)$
is a Baire space (i.e. the intersection countably many dense open sets in it is dense). In
particular, the $\tau^f_\infty$-topology on $S(A)$ is Baire.
\end{remark}

\begin {remark}\label{final_remark}
If we could show that for all $a,A\subseteq B\subseteq C$, $$\snonfork{a}{C}{A}\Rightarrow
\snonfork{a}{C}{B},$$ then this would imply that for $A\subseteq B$, $\snonfork{a}{B}{A}$ implies
$SU_s(a/A)=SU_s(a/B)$. Thus by Remark \ref{baire remark}, a Baire categoricity argument applying
Theorem \ref{uni_pcft}, will imply the existence of a bounded finite $SU_s$-rank unbounded
$\tau^f_\infty$-open set in any countable imaginary unidimensional simple theory and thus
supersimplicity will follow by Lemma \ref{main lemma}. Unfortunately, this seems to be false for a
general simple theory without stable forking.
\end {remark}

\section{$\tilde\tau^f$ and $\tilde\tau^f_{st}$-sets}
The problem of finding an unbounded $\tau^f_\infty$-open set of bounded finite $SU_s$-rank in a
countable imaginary simple unidimensional theory looked simple at first. Indeed, a Baire
categoricity argument using the "independence relation" $\snonforkempty{}{}$, instead of
$\nonforkempty{}{}$ seemed very natural but, as indicated in Remark \ref{final_remark}, doesn't
seem to work. The attempt to find other "independence relation" that is weaker than the usual
independence relation, sufficiently definable, and preserving the $SU_s$-rank seemed very
problematic too. The resolution of this obtained by analyzing sets of the form $U_{f,n}=\{a\in
\CC^s \vert\ SU_{se}(f(a))\geq n\}$, where $n<\omega$ and $f$ is an $\emptyset$-definable function
($SU_{se}$ is a variation of $SU_s$ and will be defined later). The complements of these sets
appears naturally as we assume unidimensionality; indeed, for every $a\in\CC\backslash
acl(\emptyset)$ there exists $a'\in dcl^{eq}(a)\backslash acl^{eq}(\emptyset)$ such that
$SU(a')<\omega$ and in particular $SU_{se}(a')<\omega$. The sets we will analyze in this section,
called $\tilde\tau^f$-sets, are generalizations of local versions of the sets $U_{f,n}$. The
theorem which will be crucial for the main result is that in a simple theory in which the extension
property is first-order, any minimal unbounded fiber of a $\tilde\tau^f$- set is a $\tau^f$-open
set. In this section $T$ is assumed to be a simple theory. We work in $\CC$.

\begin{definition}\em
A relation $V(x,z_1,...z_l)$ is said to be a \em pre-$\tilde\tau^f$-set relation \em if there are
$\theta(x,\tilde x,z_1,z_2,...,z_l)\in L$ and $\phi_i(\tilde x,y_i)\in L$ for $0\leq i\leq l$ such
that for all $a,d_1,...,d_l\in \CC$ we have
$$V(a,d_1,...,d_l)\ \mbox{iff}\ \exists \tilde a\
 [\theta(a,\tilde a, d_1,d_2,...,d_l)\wedge\bigwedge^{l}_{i=0} (\phi_i(\tilde a,y_i)\ \mbox{forks\
over}\ d_1d_2...d_i)]\ $$ (for $i=0$ the sequence $d_1d_2...d_i$ is interpreted as $\emptyset$). If
each $\phi_i(\tilde x,y_i)$ is assumed to be stable, $V(x,z_1,...z_l)$ is said to be a \em
pre-$\tilde\tau_{st}^f$-set relation.
\end{definition}

\begin{definition}\em
1) A \em $\tilde\tau^f$-set (over $\emptyset$) \em is a set of the form
$$\UU=\{a \vert\ \exists d_1,d_2,...d_l\
V(a,d_1,...,d_l)\}$$ for some \em pre-$\tilde\tau^f$-set relation $V(x,z_1,...z_l)$\em.


\noindent 2) A $\tilde\tau^f_{st}$-set is defined in the same way as a $\tilde\tau^f$-set but we
add the requirement that $V(x,z_1,...z_l)$ is a \em pre-$\tilde\tau_{st}^f$-set relation.\em
\end{definition}

We will say that the formula $\phi(x,y)\in L$ is \em low in $x$ \em if there exists $k<\omega$ such
that for every $\emptyset$-indiscernible sequence $(b_i \vert i<\omega)$, the set $\{\phi(x,b_i)
\vert i<\omega\}$ is inconsistent iff every subset of it of size $k$ is inconsistent. Note that
every stable formula $\phi(x,y)$ is low in both $x$ and $y$.

\begin{remark}\label{low type def}
Note that if $\phi(x,y)\in L$ is low in $x$ then the relation $F_\phi$ defined by $F_\phi(b,A)$ iff
$\phi(x,b)$ forks over $A$ is type-definable. Thus every pre-$\tilde\tau_{st}^f$-set relation is
type-definable and every $\tilde\tau_{st}^f$-set is type-definable.
\end{remark}

\begin{lemma}\label{tilde-tau-lemma}
Assume the extension property is first-order in $T$. Let $\theta(x,z_1,...,z_n)$ be a Stone-open
relation over $\emptyset$ and let $\phi_j(x,y_j)\in L$ for $j=0,..,n$. Let $U$ be the following
invariant set. For all $d_1\in \CC$, $U(d_1)$  iff

$$\exists a\exists d_2...d_n[\theta(a,d_1,...d_n)\wedge\bigwedge_{j=0}^n \phi_j(a,y_j) {\ forks\ over\
d_1...d_j}].$$ Then $U$ is a $\tau^f$-open set over $\emptyset$. If each $\phi_j(x,y_j)$ is assumed
to be low in $y_j$ and $\theta$ is assumed to be definable, then $U$ is a basic
$\tau^f_{\infty}$-open set.
\end{lemma}

\proof We prove the lemma by induction on $n\geq 1$. Consider the negation $\Gamma$ of $U$:
$$\Gamma(d_1)\ \mbox{iff}\ \forall a\forall d_2...d_n(\theta(a,d_1,...d_n)\rightarrow \bigvee_{j=0}^n
\phi_j(a,y_j)\ \mbox{dnfo}\ d_1...d_j)$$ (where "dnfo"=doesn't fork over).

\begin{subclaim}\label{Q generalized}
Let $\Gamma'$ be defined by $\Gamma'(d_1)$ iff $$\bigwedge_{\{\eta_j\}_{j=0}^{n-1}\in L}\forall
d_2...d_n[(\bigwedge_{j=0}^{n-1} \eta_j(d_1...d_n,y_j)\ \mbox{forks\ over}\ d_1...d_j)\rightarrow
\forall a \Lambda(a,d_1,...,d_n)].$$ where $\Lambda$ is defined by
$$\Lambda(a,d_1,...d_n)\ \ \mbox{iff}\ \ \theta(a,d_1,...d_n)\rightarrow \bigvee_{j=0}^n \phi_j(a,y_j)\wedge
\neg\eta_j(d_1...d_n,y_j)\ \mbox{dnfo}\ d_1...d_n$$ where $\eta_n$ denotes a contradiction. Then
$\Gamma'=\Gamma$.
\end{subclaim}

\proof Assume $\Gamma(d_1)$. Let $\eta_0,...\eta_{n-1}\in L$ and let $d_2,...d_n\in \CC$. Assume
$\eta_j(d_1...d_n,y_j)$ forks over $d_1...d_j$ for all $0\leq j\leq n-1$, and let $a\in\CC$ be such
that $\theta(a,d_1,...d_n)$.
By the assumption, we may assume $\phi_{j_0}(a,y_{j_0})$ doesn't fork over $d_1...d_{j_0}$ for some
$0\leq j_0\leq n-1$. Let $c_{j_0}$ be such that $\phi_{j_0}(a,c_{j_0})$ and
$\nonfork{a}{c_{j_0}}{d_1...d_{j_0}}$. By extension we may assume
$\nonfork{ad_1...d_n}{c_{j_0}}{d_1...d_{j_0}}$. Since $\eta_{j_0}(d_1...d_n,y_{j_0})$ forks over
$d_1...d_{j_0}$, we know $\neg\eta_{j_0}(d_1...d_n,c_{j_0})$. Therefore
$\phi_{j_0}(a,y_{j_0})\wedge \neg\eta_{j_0}(d_1...d_n,y_{j_0})$ doesn't fork over $d_1...d_{j_0}$
and in particular doesn't fork over $d_1...d_n$. Assume now $\Gamma'(d_1)$. Let $a,d_2,...d_n\in
\CC$ and assume $\theta(a,d_1,...d_n)$. It is sufficient to show that for all $0\leq j\leq n-1$ if
$\phi_j(a,y_j)$ forks over $d_1...d_j$, then there exists $\eta_j$ such that
$\eta_j(d_1...d_n,y_j)$ forks over $d_1...d_j$ and $\phi_j(a,y_j)\wedge \neg\eta_j(d_1...d_n,y_j)$
forks over $d_1....d_n$. Assume otherwise. Fix $j$, so $\phi_j(a,y_j)$ forks over $d_1...d_j$ and
$\phi_j(a,y_j)\wedge \neg\eta_j(d_1...d_n,y_j)$ doesn't fork over $d_1....d_n$ for all $\eta_j$
such that $\eta_j(d_1...d_n,y_j)$ forks over $d_1...d_j$. Let
$$\Psi(y_j)\equiv\bigwedge_{\eta_j\in F_j, \mu_j\in E_j}
\phi_j(a,y_j)\wedge\neg\eta_j(d_1...d_n,y_j)\wedge\neg\mu_j(ad_1...d_n,y_j)$$ where $$F_j=\{\eta_j
\vert\ \eta_j(d_1...d_n,y_j)\ \mbox{forks\ over}\ d_1...d_j\},$$ and $$E_j=\{\mu_j\vert\
\mu_j(ad_1...d_n,y_j)\ \mbox{forks\ over}\ d_1...d_n\}.$$ By our assumption and compactness,
$\Psi(y_j)$ is consistent. Let $c_j\models\Psi(y_j)$. Then $\phi_j(a,c_j)$,
$\nonfork{d_1...d_n}{c_j}{d_1...d_j}$, and $\nonfork{ad_1...d_n}{c_j}{d_1...d_n}$. By transitivity,
$\nonfork{ad_1...d_n}{c_j}{d_1...d_j}$. A contradiction to the assumption that $\phi_j(a,y_j)$
forks over $d_1...d_j$. The proof of Subclaim \ref{Q generalized} is complete.\\

\noindent Since the extension property is first-order for $T$, the relation $\Lambda_0$ defined by
$\Lambda_0(d_1,...d_n)\equiv\forall a \Lambda(a,d_1,...d_n)$ is type-definable. Now, clearly for
all $d_1$, $\Gamma'(d_1)$ iff $$\bigwedge_{\{\eta_j\}_{j=0}^{n-1}\in L}\forall
d_2...d_n(\neg\Lambda_0(d_1,...,d_n)\rightarrow \bigvee_{j=0}^{n-1} \eta_j(d_1...d_n,y_j)\
\mbox{dnfo}\ d_1...d_j).$$ Now, if $n=1$ then this is clearly $\tau^f$-closed. If $n>1$, then we
finish by the induction hypothesis.

\begin{corollary}\label{tau_cor}
Assume the extension property is first-order in $T$. Let $m\leq l<\omega$ and let
$d^*_1,...d^*_m\in\CC$. Let $\theta\in L$ and $\phi_i\in L$ for $i\leq l$. Let $V$ be defined by
$$V(a,d_1,...,d_l)\ \mbox{iff}\ [\theta(a,d_1,d_2,...,d_l)\wedge\bigwedge^{l}_{i=0} (\phi_i(a,y_i)\
\mbox{forks\ over}\ d_1d_2...d_i)].$$ Then the set $U$ defined by
$$U(d_{m+1})\ \mbox{iff}\ \exists a\exists d_{m+2}...d_l\ V(a,d_1^*,...d^*_m,d_{m+1},...d_l)$$ is a
$\tau^f$-open set over $d^*_1...d^*_m$.
\end{corollary}

\proof By Fact \ref{tau extensions}, there are formulas $\{\psi_j(\tilde x, w_j)\in
L(d^*_1...d^*_m)\}_{j\in J}$ such that
$$\forall a\ [\bigwedge^{m}_{i=0} (\phi_i(a,y_i)\ \mbox{forks\ over}\
d^*_1d^*_2...d^*_i)\ \mbox{iff} \bigvee_{j\in J} (\psi_j(a, w_j) \ \mbox{forks\ over}\
d^*_1d^*_2...d^*_m)].$$ Therefore by Lemma \ref{tilde-tau-lemma} (since by Lemma \ref{ext over
sets}, the extension property is first-order over $\bar d^*$ as well) $U$ is a union over $j\in J$
of $\tau^f$-open sets over $d^*_1d^*_2...d^*_m$.

\begin{theorem}\label {tilde top thm}
Assume the extension property is first-order in $T$. Then

\noindent 1) Let $\UU$ be an unbounded $\tilde\tau^f$-set over $\emptyset$. Then there exists an
unbounded $\tau^f$-open set $\UU^*$ over some finite set $A^*$ such that $\UU^*\subseteq \UU$. In
fact, if $V(x,z_1,...,z_l)$ is a pre-$\tilde\tau^f$-set relation such that $\UU=\{a\vert \exists
d_1...d_l V(a,d_1,...,d_l)\}$, and $(d^*_1,...,d^*_m)$ is any maximal sequence (with respect to
extension) such that $\exists d_{m+1}...d_l V(\CC,d^*_1,...,d^*_m,d_{m+1},...,d_l)$ is unbounded,
then
$$\UU^*=\exists d_{m+1}...d_l V(\CC,d^*_1,...,d^*_m,d_{m+1},...,d_l)$$ is a $\tau^f$-open set over
$d^*_1...d^*_m$.

\noindent 2) Let $\UU$ be an unbounded $\tilde\tau_{st}^f$-set over $\emptyset$. Then there exists
an unbounded $\tau_\infty^f$-open set $\UU^*$ over some finite set $A^*$ such that $\UU^*\subseteq
\UU$. In fact, if $V(x,z_1,...,z_l)$ is a pre-$\tilde\tau^f_{st}$-set relation such that
$\UU=\{a\vert \exists d_1...d_l V(a,d_1,...,d_l)\}$, and $(d^*_1,...,d^*_m)$ is any maximal
sequence (with respect to extension) such that $\exists d_{m+1}...d_l
V(\CC,d^*_1,...,d^*_m,d_{m+1},...,d_l)$ is non-algebraic, then $$\UU^*=\exists d_{m+1}...d_l
V(\CC,d^*_1,...,d^*_m,d_{m+1},...,d_l)$$ is a basic $\tau^f_\infty$-open set over $d^*_1...d^*_m$.
\end{theorem}

\proof By Remark \ref{low type def}, (2) is an immediate corollary of (1). It suffices, of course,
to prove the second part of (1). $T$ is PCFT by Corollary \ref{ext pcft}. Let $\bar
d^*=d^*_1...d^*_m$. First, if $m=l$ then the assertion follows immediately by Fact \ref {tau
extensions}. So, we may assume $m<l$. By maximality of $\bar d^*$, we know $\exists d_{m+2}...d_l
V(\CC,d^*_1,...,d^*_m,d'_{m+1},d_{m+2},...d_l)$ is bounded (equivalently, a union of algebraic sets
over $\bar d^*$) for every $d'_{m+1}$. Thus for every $a\in \UU^*$, there exist $\chi_a(x,\bar
z^*,z)\in L$, $k=k(\chi_a)<\omega$ and $d'_{m+1}(a)\in\CC$, such that $\forall z\forall \bar
z^*\exists^{=k} x\chi_a(x,\bar z^*,z)$ (*1) and $V(a,d^*_1,...,d^*_m,d'_{m+1}(a),d_{m+2},...d_l)$
for some $d_{m+2},...d_l\in\CC$ and $\chi_a(x,\bar d^*,d'_{m+1}(a))$ isolates the type $tp(a/\bar
d^*,d'_{m+1}(a))$.  Let $\Xi=\{\chi_a\}_{a\in \UU^*}$. For $\chi\in\Xi$, let $k=k(\chi)$ and let
$U_\chi$ be the $\bar d^*$-invariant set defined by $U_\chi(d_{m+1})$ iff $$\exists\
\mbox{distinct}\ a_1....a_k [\bigwedge_{j=1}^k \chi(a_j,\bar d^*,d_{m+1})\wedge\bigwedge_{j=1}^k\
\exists d_{m+2}...d_l V(a_j,\bar d^*,d_{m+1},d_{m+2},...d_l)]$$

\begin{subclaim}\label{tilde tau subclaim}
$U_\chi$ is a $\tau^f$-open set over $\bar d^*$.
\end{subclaim}

\proof Let $V$ be given by: $$V(a,d_1,...,d_l)\ \mbox{iff}\ \exists \tilde a\
 [\theta(a,\tilde a, d_1,d_2,...,d_l)\wedge\bigwedge^{l}_{i=0} (\phi_i(\tilde a,y_i)\ \mbox{forks\
over}\ d_1d_2...d_i)].$$ for some $\theta,\phi_i\in L$. Since $T$ is PCFT, it is sufficient to show
that there exists a $\tau^f$-open set $W=W(x,z_{m+1},\bar d^*)$ over $\bar d^*$ such that if
$U'_\chi$ is defined by
$$U'_\chi(d_{m+1})\ \mbox{iff}\ \exists\ \mbox{distinct}\ a_1....a_k [\bigwedge_{j=1}^k \chi(a_j,\bar
d^*,d_{m+1})\wedge\bigwedge_{j=1}^k\ W(a_j,d_{m+1},\bar d^*)]$$ then $U'_\chi=U_\chi$. To show this
let $W$ be defined by: $W(a,d_{m+1},\bar d^*)$ iff $$\exists \tilde a\exists d'_{m+2}...d'_l
 [\theta(a,\tilde a, d^*_1,d^*_2,...d^*_m,d_{m+1},d'_{m+2},...d'_l)\wedge\bigwedge^{l}_{i=0}
 (\phi_i(\tilde a,y_i)\ \mbox{forks\ over}\ d'_1d'_2...d'_i)]$$ where $d'_i$ is defined in the following way:
for $1\leq i\leq m$, $d'_i$ denotes $d^*_i$, and $d'_{m+1}$ denotes $d_{m+1}a$ (and the rest are
quantified variables). First note that for all $a,d_{m+1}$ with $a\in acl(d_{m+1},\bar d^*)$,
$W(a,d_{m+1},\bar d^*)$ iff $\exists d_{m+2}...d_l V(a,\bar d^*,d_{m+1},d_{m+2},...d_l)$. Thus by
(*1), $U'_\chi=U_\chi$. By Corollary \ref{tau_cor}, $W$ is a $\tau^f$-open set over $\bar d^*$. So,
the proof of Subclaim \ref{tilde tau subclaim} is complete.\\

\noindent Now, for each $\chi\in \Xi$ define $Y_\chi(x)\equiv\exists d_{m+1} (\chi(x,\bar
d^*,d_{m+1})\wedge U_\chi(d_{m+1}))$. Since $T$ is PCFT, Subclaim \ref{tilde tau subclaim} implies
$Y_\chi$ is a $\tau^f$-open set over $\bar d^*$. Note that by the definition of $U_\chi$ and (*1),
$Y_\chi\subseteq\UU^*$ for all $\chi\in \Xi$. Now, if $a\in \UU^*$, then by the choice of
$d'_{m+1}(a)$, $\chi_a$ and $k=k(\chi_a)$, we have $\chi_a(a,\bar d^*,d'_{m+1}(a))\wedge
U_{\chi_a}(d'_{m+1}(a))$. Thus $a\in Y_{\chi_a}$. Hence $\UU^*=\bigcup_{\chi\in \Xi} Y_\chi$, and
so $\UU^*$ is a $\tau^f$-open set over $\bar d^*$. The proof of Theorem \ref{tilde top thm} is
complete.

\section{Main Result}
We apply the theorem in section 8 to prove a new theorem for countable simple theories in which the
extension property is first-order. The theorem says the assumption that every non-algebraic element
has a non-algebraic element of finite $SU_{se}$-rank (a variation of the $SU_s$-rank) in its
definable closure implies the existence of an unbounded $\tau^f_\infty$-open set of bounded finite
$SU_{se}$-rank. It is here that we apply compactness, indeed this is possible because we require
our set to be only of bounded finite $SU_{se}$-rank rather than of bounded finite $SU$-rank. By the
reduction in section 7 and a corollary of section 3 the existence of such a set implies the main
result. In this section $T$ is assumed be a simple theory and we work in $\CC$ unless otherwise
stated.

\begin{remark}\label{uni_eq}
Note that by passing from $\CC$ to $\CC^{eq}$ (and vise versa) simplicity, supersimplicity and
unidimensionality are preserved (unidimensionality is less trivial, see \em [Claim 5.2, S1]\em ).
\end{remark}

\begin{definition}
1) For $a\in \CC$ and $A\subseteq \CC$ the $SU_{se}$-rank is defined by induction on $\alpha$: if
$\alpha=\beta+1$, $SU_{se}(a/A)\geq \alpha$ if there exist $B_1\supseteq B_0\supseteq A$ such that
$\sfork{a}{B_1}{B_0}$ and $SU_{se}(a/B_1)\geq\beta$. For limit $\alpha$, $SU_{se}(a/A)\geq\alpha$
if $SU_{se}(a/A)\geq\beta$ for all $\beta<\alpha$.

\noindent 2) Let $\UU$ be an $A$-invariant set. We write $SU_{se}(\UU)=\alpha$ (the $SU_{se}$-rank
of $\UU$ is $\alpha$) if $Max\{SU_{se}(p) \vert p\in S(A), p^\CC\subseteq\UU\}=\alpha$. We say that
$\UU$ has bounded finite $SU_{se}$-rank if for some $n<\omega$, $SU_{se}(\UU)=n$.
\end{definition}

\begin{remark} \label {SU_s SU_se}
Note that $SU_{se}(a/B)\leq SU_{se}(a/A)$ for all $a\in \CC$ and $A\subseteq B\subseteq\CC$ (this
is the reason for introducing $SU_{se})$. Also, clearly $SU_s(a/A)\leq SU_{se}(a/A)\leq SU(a/A)$
for all $a,A$. Clearly $SU_{se}(a/A)=0$ iff $SU_s(a/A)=0$ iff $a\in acl(A)$ for all $a,A$.
\end{remark}


\begin{theorem}\label{model baire}
Let $T$ be a countable simple theory in which the extension property is first-order and assume
$Lstp=stp$ over sets. Let $s$ be a sort such that $\CC^s$ is not algebraic. Assume for every
$a\in\CC^s\backslash acl(\emptyset)$ there exists $a'\in dcl(a)\backslash acl(\emptyset)$ such that
$SU_{se}(a')<\omega$. Then there exists an unbounded $\tau_{\infty}^f$-open set $\UU$ over a finite
set such that $\UU$ has bounded finite $SU_{se}$-rank.
\end{theorem}

\proof By a way of contradiction assume the non-existence of an unbounded $\tau_{\infty}^f$-open
set of bounded finite $SU_{se}$-rank over a finite set. It will be sufficient to show $\exists
a^*\in \CC^s\backslash acl(\emptyset)$ such that for every $\emptyset$-definable function $f$,
either $f(a^*)\in acl(\emptyset)$ or $SU_{se}(f(a^*))\geq \omega$. To show this, for every
$\emptyset$-definable function $f$ and $n<\omega$, let
$$S_{f,n}=\{a\in \CC^s \vert\ 0<SU_{se}(f(a))<n\}.$$

\begin{subclaim}\label {main}
For every non-empty $\tilde\tau_{st}^f$-set $\UU\subseteq\CC^s$ (with $\UU\cap
acl(\emptyset)=\emptyset)$ for all $\emptyset$-definable function $f$, and $n<\omega$, there exists
a non-empty $\tilde\tau_{st}^f$-set $\UU^*\subseteq \UU\cap (\CC^s\backslash S_{f,n})$.
\end{subclaim}

\noindent Assuming Subclaim \ref {main} is true, let $((f_i, n_i) \vert i<\omega)$ be an
enumeration of all such pairs $(f,n)$. By induction, let $\UU_0=\CC^s\backslash acl(\emptyset)$,
and let $\UU_{i+1}\subseteq \UU_i\cap (\CC\backslash S_{f_i,n_i})$ be a non-empty
$\tilde\tau_{st}^f$-set. Since each $\UU_i$ is type-definable, by compactness $\bigcap_{i<\omega}
\UU_i\neq \emptyset$. So, any $a^*\in \bigcap_{i<\omega} \UU_i$ will work.\\

\noindent \textbf{Proof of Subclaim \ref{main}:} Let $\UU$, $(f,n)$ be as in Subclaim \ref{main}.
Now, if there exists $a\in\UU$ such that $f(a)\in acl(\emptyset)$, let $\chi(x)\in L$ be algebraic
such that $\chi(f(a))$. By letting $\UU^*=\{a\in\UU \vert\ \models\chi(f(a))\}$ we are done. Hence
we may assume $f(a)\not\in acl(\emptyset)$ for every $a\in\UU$. Let $V(x,z_1,...z_l)$ be a
pre-$\tilde\tau^f_{st}$-set relation such that
$$\UU=\{a \vert\ \exists d_1,d_2,...d_l\ V(a,d_1,...,d_l)\}.$$ where $V$ is defined by:
$$V(a,d_1,...,d_l)\ \mbox{iff}\ \exists \tilde a\
 [\theta(a,\tilde a, d_1,d_2,...,d_l)\wedge\bigwedge^{l}_{i=0} (\phi_i(\tilde a,y_i)\ \mbox{forks\
over}\ d_1d_2...d_i)]$$ for some $\theta(x,\tilde x,z_1,z_2,...,z_l)\in L$ and stable
$\phi_i(\tilde x,y_i)\in L$ for $0\leq i\leq l$. Now, let $V_f$ be defined by: for all
$b,d_1,...,d_l\in \CC$, $$V_f(b,d_1,...,d_l)\ \mbox{iff}\ \exists a (b=f(a)\wedge
V(a,d_1,...,d_l)).$$ Then, clearly $V_f$ is a pre-$\tilde\tau_{st}^f$-set relation.
Let $\bar d^*=(d_1^*,...,d_m^*)$ be a maximal sequence, with respect to extension, ($m\leq l$) such
that $$\tilde V_f(v)\equiv \exists d_{m+1},d_{m+2},...d_l V_f(v,d^*_1,...d^*_m,d_{m+1},...d_l)$$ is
non-algebraic, or equivalently unbounded (since $\UU\not =\emptyset$ and we assume $f(a)\not\in
acl(\emptyset)$ for all $a\in\UU$, the empty sequence satisfies this property). By Theorem
\ref{tilde top thm}(2), $\tilde V_f(\CC)$ is a basic $\tau_\infty^f$-open set over $\bar d^*$. By
our assumption $\tilde V_f(\CC)$ is not of bounded finite $SU_{se}$-rank. Thus there are $a^*$ and
$d^*_{m+1},...d_l^*$ such that $V(a^*,\bar d^*,d_{m+1}^*,...,d_{l}^*)$ and $SU_{se}(f(a^*)/\bar
d^*)\geq n$. Let $E=\langle (c_i^*,e_i^*) \vert 1\leq i\leq n\rangle$ be  such that
$\sfork{f(a^*)}{e^*_i}{\bar d^*c_1^*e_1^*...c^*_i}$ for all $1\leq i\leq n$ (*1). Note that since
both $dcl$ and forking have finite character, we may assume that $c_i^*,e_i^*$ are finite tuples.
Let $\tilde a^*$ be such that:
$$\theta(a^*,\tilde a^*, d^*_1,d^*_2,...,d^*_l)\wedge\bigwedge^{l}_{i=0} (\phi_i(\tilde a^*,y_i)\
\mbox{forks\ over}\ d^*_1d^*_2...d^*_i)\ \ (*2).$$ Now, by maximality of $\bar d^*$, $f(a^*)\in
acl(\bar d^*d^*_{m+1})$. By taking a non-forking extension of $tp(E/acl(\bar d^*d^*_{m+1}))$ over
$acl(d_1^*...d^*_la^*\tilde a^*)$ we may assume that $\nonfork{a^*\tilde a^*d_1^*...d_l^*}{E}{\bar
d^*d^*_{m+1}}$ and (*1) and (*2) still hold. Thus $\nonfork{a^*\tilde
a^*}{d^*_{1}...d_i^*E}{d^*_{1}...d^*_i}$ for all $m+1\leq i\leq l$. Hence by (*2), we conclude
$\phi_i(\tilde a^*,y_i)\ \mbox{forks\ over}\ d^*_1d^*_2...d^*_iE$ for all $m+1\leq i\leq l$. By
(*1) and symmetry of $\snonfork{}{}{}$ (Lemma \ref{s-sym-trans}), there are stable
$\psi_i(x_i,w_i)\in L$ and $\emptyset$-definable functions $g_i,h_i$ for $1\leq i\leq n$ such that
if $a^*_i=g_i(f(a^*),\bar d^*c_1^*e_1^*...c^*_i)$, and  $b^*_i=h_i(e^*_i,\bar
d^*c_1^*e_1^*...c^*_i)$, then $\psi_i(a^*_i,b^*_i)$ and $\psi_i(a^*_i, w_i)$ forks over $\bar
d^*c_1^*e_1^*...c^*_i$. Now, let us define a relation $V^*$ in the following way:
$$V^*(a,d_1,...d_m,c_1,e_1,..c_n,e_n,d_{m+1},..d_l)\ \mbox{iff}\ \exists \tilde a, \tilde
a'=\tilde a'_1..\tilde a'_n, \tilde b'=\tilde b'_1..\tilde b'_n (\theta^*\wedge V_0\wedge V_1\wedge
V_2)$$ where, $\theta^*$ is defined by: $\theta^*(a,\tilde a,\tilde a',\tilde
b',d_1,..d_m,c_1,e_1,..c_n,e_n,d_{m+1},..d_l)\equiv$
$$\theta(a,\tilde a,d_1,..d_l)\wedge
\bigwedge_{i=1}^n [\psi_i(\tilde a'_i,\tilde b'_i)\wedge\ (\tilde
a'_i=g_i(f(a),d_1,..d_m,c_1,e_1,..c_i)\wedge\ (\tilde b'_i=h_i(e_i,d_1,..d_m,c_1,e_1,..c_i))],$$
$V_0$ is defined by:
$$V_0(\tilde a,d_1,...d_m)\ \mbox{iff} \bigwedge^{m}_{i=0} (\phi_i(\tilde a,y_i)\ \mbox{forks\
over}\ d_1d_2...d_i),$$

\noindent $V_1$ is defined by: $$V_1(\tilde a',d_1,...d_m,c_1,e_1,...c_n,e_n)\ \mbox{iff}\
\bigwedge^{n}_{i=1} (\psi_i(\tilde a'_i,w_i)\ \mbox{forks\ over}\ d_1d_2...d_mc_1e_1...c_i),$$ and
$V_2$ is defined by: $$V_2(\tilde a,d_1,..d_m,c_1,e_1,..c_n,e_n,d_{m+1},..d_l)\
\mbox{iff}\bigwedge^{l}_{i=m+1}(\phi_i(\tilde a,y_i)\ \mbox{forks\ over}\
d_1d_2..d_mc_1e_1..c_ne_nd_{m+1}..d_i).$$ Note that $V^*$ is a pre-$\tilde\tau_{st}^f$-set
relation. Thus $$\UU^*=\{a \vert\ \exists d_1,..d_m,c_1,e_1,..c_n,e_n,d_{m+1},..d_l\
V^*(a,d_1,..d_m,c_1,e_1,..c_n,e_n,d_{m+1},..d_l)\}$$ is a $\tilde\tau_{st}^f$-set. By the
construction of $a^*,d^*_1,..d^*_m,c^*_1,e^*_1,..c^*_n,e^*_n,d^*_{m+1},..d^*_l$,
$\UU^*\not=\emptyset$. By the definition of $\UU^*$, $\UU^*\subseteq \UU\cap(\CC^s\backslash
S_{f,n})$ (note that if $a\in\UU^*$, then there are $d_1,...,d_m\in\CC$ such that
$SU_{se}(f(a)/d_1...d_m)\geq n$ and thus by Remark \ref {SU_s SU_se}, $SU_{se}(f(a))\geq n$). So,
the proof of Subclaim \ref{main} is complete, and thus so is the proof of the theorem.

\begin{theorem}
Let $T$ be a countable imaginary simple unidimensional theory. Then $T$ is supersimple.
\end{theorem}

\proof By adding countably many constants we may assume there exists $p_0\in S(\emptyset)$ of
$SU$-rank 1 (each of the assumptions is preserved, as well as the corollary). Now, by Remark \ref
{uni_eq} we may work in $\CC^{eq}$. Fix a non-algebraic sort $s$. Since $T$ is unidimensional and
imaginary, by Fact \ref{internal} for every $a\in \CC^s\backslash acl(\emptyset)$ there exists
$a'\in dcl(a)\backslash acl(\emptyset)$ such that $tp(a')$ is $p_0$-internal; thus $SU(a')<\omega$
and in particular $SU_{se}(a')<\omega$. By Corollary \ref{uni ext prop} the extension property is
first-order in $T$. By Theorem \ref{model baire}, there exists an unbounded $\tau_{\infty}^f$-open
set $\UU$ over a finite set such that $\UU$ has bounded finite
$SU_{se}$-rank. By Lemma \ref {main lemma}, $T$ is supersimple.\\

Recall that a theory $T$ has the \em wnfcp(=weak non finite cover property) \em if for each
$L$-formula $\phi(x,y)$, the $D_{\phi}$-rank is finite (equivalently, $\phi(x,y)$ is low in $x$)
and definable (the $D_{\phi}$-rank of a formula $\psi(x,a)$ is defined by: $D_{\phi}(\psi(x,a))
\geq 0$ if $\psi(x,a)$ is consistent; $D_{\phi}(\psi(x,a)) \geq \alpha +1$ if for some $b$,
$D_{\phi}(\psi(x,a)\wedge\phi(x,b))\geq \alpha$ and $\phi(x,b)$ divides over $a$; and for limit
$\delta$, $D_{\phi}(\psi(x,a))\geq \delta$ if it is $\geq \alpha$ for all $\alpha<\delta$).

\begin{corollary}
Let $T$ be a countable imaginary simple unidimensional theory. Then $T$ is low and thus has the
wnfcp.
\end{corollary}

\proof By Fact \ref{supersimple definable}, $T$ has bounded finite $SU$-rank in any given sort.
Thus the global $D$-rank of any sort is finite. Now, let $\phi(x,y)\in L$. Then $\phi(x,y)$ is low
in $x$ iff $Sup\{D(x=x,\phi(x,y),k) \vert\ k<\omega\}<\omega$. So, clearly every $\phi(x,y)$ is low
in $x$. Thus $T$ is low. By Corollary \ref{uni ext prop} the extension property is first-order in
any unidimensional theory. We conclude $T$ has the wnfcp (see [BPV], Corollary 4.6).

\end{document}